\documentclass[AMA,Times1COL]{WileyNJDv5} 

\usepackage{amssymb,amsmath,mathtools}
\usepackage{amsthm}
\usepackage{stmaryrd}
\usepackage{setspace}
\articletype{Article Type}%

\received{Date Month Year}
\revised{Date Month Year}
\accepted{Date Month Year}
\journal{Journal}
\volume{00}
\copyyear{2024}
\startpage{1}

\begin{document}

\title{A deep learning-based surrogate model for seismic data assimilation in fault activation modeling}

\author[1]{Caterina Millevoi}

\author[1]{Claudia Zoccarato}

\author[1]{Massimiliano Ferronato}

\authormark{MILLEVOI \textsc{et al.}}
\titlemark{SURMODEL}

\address[1]{\orgdiv{Department of Civil, Environmental and Architectural Engineering}, \orgname{University of Padova}, \orgaddress{\state{Padova}, \country{Italy}}}

\corres{Corresponding author Caterina Millevoi, 
\email{caterina.millevoi@unipd.it}}

\presentaddress{Department of Civil, Environmental and Architectural Engineering, University of Padova, via Marzolo 9, 35131 Padova (PD), Italy.}

\abstract[Abstract]{Assessing the safety and environmental impacts of subsurface resource exploitation and management is critical and requires robust geomechanical modeling. However, uncertainties stemming from model assumptions, intrinsic variability of governing parameters, and data errors challenge the reliability of predictions. In the absence of direct measurements, inverse modeling and stochastic data assimilation methods can offer reliable solutions, but in complex and large-scale settings, the computational expense can become prohibitive.

To address these challenges, this paper presents a deep learning-based surrogate model (SurMoDeL) designed for seismic data assimilation in fault activation modeling. The surrogate model leverages neural networks to provide simplified yet accurate representations of complex geophysical systems, enabling faster simulations and analyses essential for uncertainty quantification. The work proposes two different methods to integrate an understanding of fault behavior into the model, thereby enhancing the accuracy of its predictions. The application of the proxy model to integrate seismic data through effective data assimilation techniques efficiently constrains the uncertain parameters, thus bridging the gap between theoretical models and real-world observations.
}


\maketitle

\renewcommand\thefootnote{}
\footnotetext{\textbf{Abbreviations:} MCMC, Markov Chain Monte Carlo; DL, Deep Learning; NN, Neural Network; QoI, Quantity of Interest; KKT, Karush-Kuhn-Tucker; FE, Finite Element; SGD, Stochastic Gradient Descent; pdf, posterior distribution function.
}

\renewcommand\thefootnote{\fnsymbol{footnote}}
\setcounter{footnote}{1}

\section{Introduction}\label{sec:intro}

Surrogate models have become important tools in several applications, especially in multi-scale and multi-physics scenarios involving high uncertainties and complex simulations. In essence, these models can provide simplified representations of complex systems, enabling faster simulations and analyses especially when ensembles of realizations are needed for the sake of uncertainty quantification purposes. 

In the context of geomechanical subsurface simulations, surrogate models 
have been employed to investigate the poroelasticity problem with random coefficients~\cite{Bot_etal19}, predict and quantify the uncertainty of land subsidence  models~\cite{BotRos2017,Zoc_etal20,Gaz_etal21,Gaz_etal23}, analyze the sensitivity factors controlling earth fissures due to overexploitation of groundwater resources~\cite{Li_etal22}, 
approximate the contact mechanics problem~\cite{CasJhaJua16,Zoc_etal19}, and perform global sensitivity analysis in geomechanical fractured reservoirs and hydraulically fractured wells~\cite{Verde_2015,RezNakSid_etal_2020}.
Among all cited case studies, the presence of faults within the geological formations introduces significant challenges. These challenges arise from the discontinuous nature of the problem and the complex interactions between mechanical and hydraulic processes. This leads to high uncertainty, for example, in the reservoir geology, the pore-pressure distribution, and the fault hydro-mechanical properties~\cite{CasJhaJua16,Bla_etal20}. 

Fault activation and generation of fractures are caused by stress changes due to injection and/or production of fluids into and from the surface. This activity could affect the reservoir formation integrity and cause several environmental hazards, such as fluid leakage, land motion, and induced seismic events~\cite{Rut_etal13,Rut_etal16,Bla_etal22}. Therefore, the generation and use of reliable models to forecast and prevent injection-induced fault motion and the consequently triggered seismicity with possible permanent damage is of utmost importance.

For this reason, the process of data assimilation, which involves integrating observational data into models to improve their accuracy and reliability, can be an important tool in fault modeling for updating the model parameters and the model states based on the latest available data. To this end, effective data assimilation can help bridge the gap between theoretical models and real-world observations, enhancing the model ability to forecast fault activation and the associated seismic risks. This, in turn, aids in better risk management and decision-making in the context of subsurface resource exploitation.

There are several methods for assimilating seismic and geophysical data into geomechanical models. Chang et al.~\cite{ChaCheZha10} used the ensemble Kalman filter to estimate reservoir flow and material properties by jointly assimilating dynamic flow and geomechanical observations. Emerick and Reynolds proposed a multiple assimilation of time-lapse seismic data to improve the ensemble Kalman Filter~\cite{EmeRey12} and used the ensemble smoother multiple data assimilation to generate multiple realizations of the porosity, net-to-gross ratio and permeability fields by history matching production and seismic impedance data~\cite{EmeRey13}. 
Luo et al. implemented a wavelet-based sparse representation procedure for 2D~\cite{Luo_etal16} and 3D~\cite{Luo_etal18} seismic data assimilation problems. Nejadi et al.~\cite{Nej_etal19} incorporated data matching at the well locations in a Bayesian inversion framework and constrained the model space by using a seismic impedance volume to estimate physically plausible porosity distributions with ensemble-based Markov Chain Monte Carlo (MCMC) approach.
The majority of these methods need repetitive forward simulations to generate prior ensembles of realizations, which can be unfeasible in terms of the computational cost for large scale and complex systems.
The need of fast and reliable predictions is therefore critical in ensemble-based data assimilation techniques. However, recently implemented techniques such as polynomial-based proxy models can struggle to accurately capture the behavior of faults. Indeed, the discontinuous processes associated with fault activation, such as sudden slips and changes in permeability, are particularly difficult to model~\cite{Zoc_etal19}. 

This limitation requires the development of different surrogate modeling techniques capable of handling such complexities.
Deep learning (DL)-based surrogate models have shown significant promise in the field of porous media~\cite{Pan_etal14,She_etal22, Qi_etal23}. By leveraging large datasets and powerful neural network (NN) architectures, DL models can learn complex patterns and relationships within the data. This capability makes them well-suited for modeling also the intricate dynamics of fault activation in poromechanics~\cite{Meg_etal23}. In fact, data-driven approaches such as NNs and other machine learning algorithms can be trained on seismic and geophysical data to develop predictive proxy models for fault activation and can then be integrated with traditional geomechanical models for enhanced predictions~\cite{Lu_etal23}.

In this paper, we propose a novel DL-based surrogate model (SurMoDeL) specifically designed for data assimilation in fault activation modeling. Our model is trained on a realistic dataset, simulating a discontinuous process that includes fault opening events due to excessive groundwater pumping. One of the key innovations of our approach is its ability to handle discontinuities effectively, since the DL-based model incorporates a physics-informed mechanism that makes it aware of the fault behavior. The proposed method is capable of detecting how probable is the occurrence of fault opening and integrate this information in building the surrogate solution. The use of this physical principle into the DL model ensures more accurate and reliable predictions.
Moreover, the use of a Bayesian-based MCMC method combined with the proposed surrogate model and seismic data assimilation, offers an efficient approach to parameter estimation in complex geomechanical models. 

The application to the 3D synthetic test case demonstrates the method ability to update model parameters using seismic data, highlighting the importance of data for uncertainty reduction and the effectiveness of the SurMoDeL in mimicking the outcome of the full order model and reducing the computational demand. This development can potentially improve our understanding and prediction of geological processes, leading to better management and mitigation of risks associated with fault activation.

The paper is organized as follows:
Section~\ref{sec:geom} describes the implemented workflow, including the full forward model of fault activation and its surrogate approximation by DL. The application set up to a 3D synthetic test case where fluid is pumped from a 1100-m-deep faulted reservoir is presented in Section~\ref{sec:setup}. Training and validation of the SurMoDeL are discussed in Section~\ref{sec:training}. Section~\ref{sec:faultopen} introduces two approaches to incorporate fault behavior awareness into the proxy model. Section~\ref{sec:SA} outlines a global sensitivity analysis of the model inputs, whereas the Bayesian inversion results for parameter estimations are described in Section~\ref{sec:DA}. A closing section concludes the paper.

\section{Fault activation modeling}\label{sec:geom}

Fault activation is a critical issue in the context of subsurface resource management, such as hydrocarbon extraction or storage, but also geothermal energy and groundwater production. When the stress state within a geological formation exceeds a certain failure criterion, pre-existing faults can become active, leading to potentially hazardous slip events and high energy dissipation. Therefore, predicting fault activation is essential for mitigating risks associated with induced seismicity.

The numerical simulation of fault activation, as well as many other subsurface processes, is subject to a number of uncertainties. For example, knowledge of the geometry and heterogeneity of deep formations is crucial to obtain a reliable modeling result. Similarly, material parameters and governing constitutive laws are often very uncertain and can lead to a broad range of possible outcomes, especially in a strongly non-linear model. This is particularly true when faults are involved because their physical characterization can hardly be supported by direct measurements.
Seismic data can play an important role to reduce the uncertainties connected to the fault characterization. In this regard, seismic monitoring networks, which allow localizing the events and quantifying the seismic moment, can provide insights into the subsurface stress state and fault mechanics, offering real-time or near-real-time observations of micro-seismic events before a potentially big occurrence. This information is crucial for understanding the conditions under which faults might slip and for developing geomechanical models that can predict such events. 

Generally speaking, we can state that the outcome $\mathbf{y}\in\mathbb{R}^K$ at every point $\mathbf{x}$ of the space domain $\Omega\subset\mathbb{R}^3$ and every instant $t$ of the time domain $[0,+\infty[$ arises from some forward model $\mathcal{S}$ providing the functional relationship between the forcing terms (loads) $\mathbf{F}$ and the independent material parameter vector $\mathbf{p}\in\mathbb{R}^n$:
\begin{equation}
    \mathbf{y}(\mathbf{x},t)=\mathcal{S}(\mathbf{F},\mathbf{p}).
    \label{eq:genfram}
\end{equation}
Seismic data represent the vector $\mathbf{q}\in\mathbf{R}^Q$ of the quantities of interest (QoIs), or observables, which are related to the model states $\mathbf{y}$ at some point of $\Omega$ by a proper mapping $\mathcal{M}:\mathbf{y}\rightarrow\mathbf{q}$, such that:
\begin{equation}
    \mathbf{q}(t)=\mathcal{M}\circ\mathcal{S}(\mathbf{F},\mathbf{p}),
    \label{eq:obs}
\end{equation}
where the parameter vector $\mathbf{p}$ is affected by some uncertainty.
Our objective is to solve the inverse problem and estimate the posterior distributions of $\mathbf{p}$ conditioned on prior knowledge and the observables $\mathbf{q}$. This can be done by using a Bayesian inference approach, where the posterior likelihood function $P(\mathbf{q}|\mathbf{F},\mathbf{p})$ is sampled by using a MCMC method.

To this aim, we need: (i) an appropriate forward model $\mathcal{S}$ to replicate the relevant physical processes, (ii) the mapping $\mathcal{M}$ that connects the outcome of the forward model with the available observables, and (iii) a fast algorithm to generate the ensembles of realizations required by the MCMC algorithm. Since the numerical simulations with the full forward model are usually very time consuming, in this work we introduce a DL-based surrogate model that can effectively replace $\mathcal{S}$.   

\subsection{Full forward model}\label{sec:numsim}

The simulation of the inception of fault activation in a geological medium is governed by frictional contact mechanics. The relative displacement between the contact surfaces can occur under particular stress conditions and evolves following specific constraints, such as the impenetrability of solid bodies and the governing static-dynamic friction law.
From a mathematical point of view, we consider the equilibrium of a deformable solid occupying the finite domain $\Omega\subset\mathbb{R}^3$ with the assumption of quasi-static conditions and infinitesimal strain. If $\Gamma_f$ denotes a pair of inner contact surfaces with normal direction $\mathbf{n}_f$, the governing linear momentum balance with the contact constraints reads \cite{KikOde88,Lau03,Wri06}:
\begin{subequations}
  \begin{align}
  & -\nabla \cdot \boldsymbol{\sigma} (\mathbf{u}) = \mathbf{b}, & & & & & & \mbox{(equilibrium)}, \label{eq:equilibrium} \\ 
  & t_N = \mathbf{t} \cdot \mathbf{n}_f \le 0,
  &
  & g_N = \llbracket \mathbf{u} \rrbracket \cdot \mathbf{n}_f \ge 0,
  &
  & t_N g_N = 0,
  &
  & \mbox{(impenetrability)},
  \label{eq:normal_contact_KKT} \\
  & \left\| \mathbf{t}_T \right\|_2 \le \tau_{\max}(t_N),
  &
  & \dot{\mathbf{g}}_T \cdot \mathbf{t}_T = \tau_{\max}(t_N) || \dot{\mathbf{g}}_T ||_2,
  &&&
  & \mbox{(friction)}.
  \label{eq:frictional_contact_KKT}
\end{align}
\label{eq:KKT}\null
\end{subequations}

In the inequality-constrained problem \eqref{eq:KKT}, the displacement $\mathbf{u}$ in $\Omega$ and the traction $\mathbf{t}$ over $\Gamma_f$ are the primary unknowns, with: $\mathbf{b}$ the external body forces; $\boldsymbol{\sigma}
(\mathbf{u})$ the stress tensor; $\mathbf{t} = t_N \mathbf{n}_f + \mathbf{t}_T$ the traction
over $\Gamma_f$ decomposed into its normal and tangential components, $t_N$ and $\mathbf{t}_T$; $\llbracket \mathbf{u} \rrbracket = 
g_N \mathbf{n}_f + \mathbf{g}_T$ the jump of $\mathbf{u}$ across $\Gamma_f$, decomposed into its normal and tangential components, $g_N$ and $\mathbf{g}_T$;
and $\tau_{\max}(t_N)$ a bounding value for the measure of $\mathbf{t}_T$. Relationships \eqref{eq:normal_contact_KKT}-\eqref{eq:frictional_contact_KKT} are the Karush-Kuhn-Tucker (KKT) complementary conditions for normal and frictional contact \cite{SimHug98}. In essence, they state that: (i) the normal traction must be compressive if the contact exists, with no penetration allowed between the two sides of the discontinuity surface $\Gamma_f$ (equation \eqref{eq:normal_contact_KKT}), and (ii) an upper bound for the magnitude of the tangential component of traction is set, at which slip is allowed and is collinear with friction (equation \eqref{eq:frictional_contact_KKT}). The mathematical problem is closed by adding the constitutive relationships for the stress $\boldsymbol{\sigma}(\mathbf{u})$ and friction $\tau_{\max}(t_N)$, and prescribing appropriate Dirichlet and Neumann boundary conditions.

In the context of the geological porous media of interest, the external body forces $\mathbf{b}$ are related to the variation of the pore pressure $p_\alpha$ for the fluid phase $\alpha$ due to human intervention. The distribution of $p_\alpha$ within $\Omega$ for every time instant $t$ is governed by the generalized multiphase flow model:
\begin{equation}
    - \nabla \cdot \left[ \frac{\kappa \rho_\alpha}{\mu_\alpha} \nabla p_\alpha \right] + \frac{\partial}{\partial t} \left(\varphi S_\alpha p_\alpha \right) = q,
    \label{eq:flow}
\end{equation}
where $\kappa$ is the permeability of the porous medium, $\rho_\alpha$ and $\mu_\alpha$ are the density and the viscosity of the fluid phase $\alpha$, $\varphi$ is the porosity and $S_\alpha$ the saturation index. Following a one-way coupled approach, the body forces $\mathbf{b}$ used in the equilibrium equation \eqref{eq:equilibrium} depend on the gradient of the equivalent pore pressure $\tilde{p}$ by the Biot coefficient $b$:
\begin{equation}
    \mathbf{b} = b \nabla \tilde{p}, \qquad \tilde{p} = \sum_\alpha S_\alpha p_\alpha.
    \label{eq:forces}
\end{equation}

A well-posed formulation of problem \eqref{eq:KKT} can be obtained by prescribing the minimization of the associated constrained variational principle in a mathematically exact way by using Lagrange multipliers \cite{Ber84,Wri06}. Convergence and numerical stability of the non linear problem is generally improved \cite{HagHueWoh08,fraferjantea16} at the cost of adding new variables as primary unknowns and increasing the overall problem size.
Lagrange multipliers have the physical meaning of traction vector $\mathbf{t}$ living on the discontinuity surface $\Gamma_f$. Denoting by $\bf{\mathcal{U}}$ and $\bf{\mathcal{U}}_0$ the subspace of $[H^1(\Omega)]^3$ acting as trial and test spaces for the displacement, respectively, and by $\bf{\mathcal{T}}(\mathbf{t})$ the appropriate function space for the Lagrange multipliers \cite{KikOde88}, the weak variational form of \eqref{eq:KKT} consists of finding $\{\mathbf{u},\mathbf{t}\}\in\bf{\mathcal{U}}\times\bf{\mathcal{T}}(\mathbf{t})$ such that:
\begin{subequations}
  \begin{align}
      \left(\nabla^s\boldsymbol{\eta},\boldsymbol{\sigma}\right)_{\Omega} + \left(\llbracket\boldsymbol{\eta}\rrbracket,\mathbf{t}\right)_{\Gamma_f} = \left(\boldsymbol{\eta},\mathbf{b}\right)_{\overline{\Omega}}, & & &\forall \boldsymbol{\eta}\in\bf{\mathcal{U}}_0, \label{eq:virtualwork} \\
      \left(t_N-\mu_N,g_N\right)_{\Gamma_f} + \left(\mathbf{t}_T-\boldsymbol{\mu}_T,\dot{\mathbf{g}}_T\right)_{\Gamma_f} \geq 0, & & &\forall \boldsymbol{\mu}\in\bf{\mathcal{T}}(\mathbf{t}), \label{eq:congruence}
  \end{align}
  \label{eq:weak_form}
\end{subequations}
where \eqref{eq:virtualwork} expresses the virtual work principle and \eqref{eq:congruence} the compatibility conditions for the contact surface. The subscripts $N$ and $T$ for the test function $\boldsymbol{\mu}$ denote the normal and tangential projection, respectively, of $\boldsymbol{\mu}$ onto $\Gamma_f$. The variational inequality \eqref{eq:congruence} can be transformed into an equality by detecting the current contact operating mode of every point lying on $\Gamma_f$, for instance with the aid of an active-set algorithm. 
According to the current operating mode, $\Gamma_f$ can be partitioned into three portions:
\begin{itemize}
    \item stick region $\Gamma_f^{stick}$: there is no discontinuity in the displacement function across the surface $\Gamma_f$ ($\llbracket\mathbf{u}\rrbracket=0$) and the traction $\mathbf{t}$ is unknown;
    \item slip region $\Gamma_f^{slip}$: the fault is stick in the normal direction ($g_N=0$ and $t_N$ is unknown), but a relative displacement between the two contact faces is allowed ($\mathbf{g}_T\neq\mathbf{0}$) with $\mathbf{t}_T=\tau_{\max}(t_N)\dot{\mathbf{g}}_T/{\left\|\dot{\mathbf{g}}_T\right\|_2}$;
    \item open region $\Gamma_f^{open}$: a free relative displacement $\llbracket\mathbf{u}\rrbracket$ is allowed with $\mathbf{t}=0$.
\end{itemize}
Dissipation of energy with the potential generation of micro-seismic events can occur only in the slip region $\Gamma_f^{slip}$, whose identification is part of the outcome of the model.

Discretization of the continuous problem \eqref{eq:weak_form} is finally carried out by replacing the mixed function space $\bf{\mathcal{U}}\times\bf{\mathcal{T}}(\mathbf{t})$ with the discrete subspace $\bf{\mathcal{U}}^h\times\bf{\mathcal{T}}^h(\mathbf{t}^h)$ associated to a conforming partition of the geometrical domain. In this work, we use a classical Finite Element (FE) discretization of the porous medium with a piecewise linear and a piecewise constant representation of $\mathbf{u}^h$ and $\mathbf{t}^h$, respectively \cite{fraferjantea16,fr2020alg,FraFerFriJan22}. 

\subsection{Parameter space and observables}\label{sec:paramspace}

The parameter space includes the uncertain quantities that influence the outcome of the geomechanical model. As already observed above, there is a variety of uncertain entities affecting the modeling results, including the formation geometry and the actual relevant physical processes. In our work, we assume that the problem geometry is deterministic and that rock deformation due to a known distribution of pore pressure variation in a deep reservoir is the main process of interest. We focus our attention mainly on the material properties that influence the outcome of the geomechanical model. For their definition, we need to introduce a constitutive relationship for both the stress tensor $\boldsymbol{\sigma}(\mathbf{u})$ and friction $\tau_{\max}(t_N)$. In terms of the former, the literature has many well-established options. Without loss of generality for the model, we limit our analysis to a standard isotropic linear elastic law defined by the values of the Young modulus and the Poisson ratio $\nu$ of the porous medium. 
By distinction, the correct reproduction of the fault dynamics is much more difficult. Actually, faults are complex three-dimensional structures consisting of a plastic core surrounded by an inner and an outer damage zone characterized by a variable distribution of rock joints. Their mechanical behavior can be idealized as that of a pair of frictional contact surfaces, as described in Section \ref{sec:numsim}, but the material parameters should take into account the general average behavior of a large rock volume and cannot be obtained by direct measurements.
For this reason, the mechanical properties that govern the friction behavior of the fault are usually much more uncertain and difficult to estimate than the material parameters of the rock constitutive law, which can be obtained by laboratory experiments and confirmed by in situ indirect measurements; see, for instance, \cite{Fer_etal13,Zoc_etal16,ZocFerTea18}. 
In our study, we focus primarily on the parameters that are central to the fault activation dynamics, but whose evaluation cannot be easily supported by direct measurements.

A well-established definition of $\tau_{\max}(t_N)$ for the prediction of fault activation is based on the classical Mohr-Coulomb failure criterion:
\begin{equation}\label{eq:mohrcoulomb}
    \tau_{\max}=\tau_0-t_N\tan\phi,
\end{equation}
where $\tau_0$ is the cohesion and $\phi$ is the friction angle. According to the contact constraints \eqref{eq:normal_contact_KKT}-\eqref{eq:frictional_contact_KKT}, when $\|\mathbf{t}_T\|_2$ reaches $\tau_{\max}$, sliding begins, and when $t_N$ goes down to 0 the fault opens.
The parameter space for the fault properties therefore includes $\tau_0$ and $\phi$.
Another important aspect controlling the possible fault activation is the initial stress regime operating on $\Gamma_f$.
Identification of the stress regime is often one of the most uncertain elements in a geomechanical simulation and is usually performed as an average over a large area because it can typically be derived from geological considerations without direct measurements.
The initial undisturbed stress tensor is defined by the principal stresses $\sigma_1$, $\sigma_2$, and $\sigma_3$, and in most situations it is acceptable to hypothesize that the undisturbed principal stress directions are almost vertical and horizontal.
Therefore, we assume for our application that the largest (in absolute value) principal stress, $\sigma_3$, is vertical, while $\sigma_1$ and $\sigma_2$ are directed towards the $x-$ and $y-$axis in a Cartesian reference frame.
Of course, should different indications be available for the specific problem at hand, they could be used in the model with no modifications to the approach presented here. 
Not all the components of the natural stress tensor are equally uncertain.
The value of the principal vertical stress $\sigma_3$ is generally characterized by a high confidence, because it can be defined as a function of the depth $z$ according to the density of the deposited sediments.
The values of the horizontal principal stresses, $\sigma_1$ and $\sigma_2$, are typically much more uncertain. 
For example, in a normally consolidated regime, the horizontal principal stress state is isotropic, with $\sigma_1$ ($=\sigma_2$) equal to a fraction of $\sigma_3$ according to the confinement factor $M=\nu/{(1-\nu)}$. 
If the stress regime is not normally consolidated, $\sigma_1\neq\sigma_2$ and we can define two different values of the confinement factor, $M_1$ and $M_2$, such that $\sigma_1=M_1\sigma_3$ and $\sigma_2=M_2\sigma_3$. An estimate of the magnitude and orientation of the minimum stress $\sigma_1$, although not as accurate as the vertical stress, can be indirectly derived from the distortion of the casing of production or monitoring wells. In contrast, a reliable estimate of the intermediate principal stress $\sigma_2$ is very difficult to obtain.
With the aim at considering an appropriate variability range for the most uncertain material properties and, at the same time, limiting the size of the parameter space, we define the set of uncertain parameters $\mathbf{p}=\{\tau_0,\phi,M_2\}$, assumed to be constant in space and time. This choice is also consistent with the same application carried out in \cite{Zoc_etal19}.

The set of observables $\mathbf{q}$ can be provided by a micro-seismic monitoring network, which measures real-time data on seismic events. Typically, these networks consist of arrays of seismometers strategically placed to detect and record the ground motion caused by an occurrence, and estimate the related energy dissipation down to a very small (even negative) magnitude. In particular, the collected data allow for computing the seismic moment, which is related to the physical properties of the fault and the slip occurring during an event. 

The seismic moment $M_0$ is a measure of the total energy released by a seismic event and is defined as \cite{KanAnd75}:
\begin{equation}\label{eq:M0}
M_0 = G \cdot A_a \cdot \delta_S,
\end{equation}
where $G$ is the shear modulus of the rock surrounding the activated fault, $A_a$ is the fault slipping area, i.e., the activated area, and $\delta_S$ is the average relative tangential displacement on the fault.
The activated area and the average fault slip are results that can be computed at each time-step $t_i$ of the simulation by means of the full forward model described above and represent the state vector $\mathbf{y}(t)=\{A_a,\delta_S\}$, while the vector of the observables is $\mathbf{q}(t)=\{M_0\}$. In essence, $A_a$ is the measure of $\Gamma_f^{slip}$ and $\delta_S$ is the integral of $\|\mathbf{g}_T\|_2$ over $A_a$: 
\begin{equation}\label{eq:output}
A_a = \left| \Gamma_f^{slip} \right|, \qquad
\delta_S = \frac{1}{A_a}\int_{\Gamma_f^{slip}} \left\| \mathbf{g}_T \right\|_2 \; d\Gamma.
\end{equation}

\subsection{Surrogate model design} 
\label{sec:DL}

In order to save computational time in the generation of ensembles of realizations with the full forward model, which is potentially very large and includes severe non-linearities, we want to design a surrogate model able to approximate the action of $\mathcal{S}$ on the loads $\mathbf{F}(t)$ and the parameters $\mathbf{p}$ to obtain the output state vector $\mathbf{y}(t)=\{A_a,\delta_S\}$ for every simulation time $t_i$. To this aim, we use basic tools in a DL framework. 
The fundamental DL unit is known as a neural network, which is a mathematical function mimicking the relationship between a set of inputs and corresponding outputs. 
This function is constructed by combining simple (nonlinear) functions, which enables the learning of complex feature hierarchies.
NNs can be used for both regression and classification tasks: in regression, the network generates continuous outputs, while in classification, it produces discrete values. In a supervised framework, the objective is to use these networks to create a model using a dataset of input-output pairs, allowing it to learn the relationship between the two and generalize to new data. This process is known as training. 
A crucial aspect in the training is the requirement of a sufficiently large amount of data, which can be obtained from measurements and investigations or specifically generated by simulations.

A feedforward NN is designed to approximate an unknown function $\mathbf{f}:\mathbb{R}^s\rightarrow\mathbb{R}^m$ using training data points. The NN approximation of $\mathbf{f}$, denoted as $\hat{\mathbf{f}}$, is achieved through the recursive composition of the function $\boldsymbol{\Sigma}^{(l)}$:
\begin{equation}\label{eq:layer}
\boldsymbol{\Sigma}^{(l)}(\mathbf{x}^{(l)}) = \sigma^{(l)}.(\mathbf{W}^{(l)}\mathbf{x}^{(l)}+\mathbf{b}^{(l)}),
\end{equation}
where $\mathbf{W}^{(l)} \in \mathbb{R}^{n_l \times n_{l-1}}$ is the matrix of the weights, $\mathbf{b}^{(l)} \in \mathbb{R}^{n_l}$ is the vector containing the biases, and $\sigma^{(l)}$ is the activation function for the $l$-th layer. The output layer is the final layer, while the preceding layers are hidden layers. The number of neurons in layer $l$ is denoted by $n_l$. Activation functions, specified by the user, typically have a limited range and are non-linear to keep the weight values low and to introduce non-linearity to the NN. The MATLAB-inspired notation $\sigma^{(l)}.(\mathbf{v})$ indicates that the function $\sigma^{(l)}$ is applied component-wise to the argument vector $\mathbf{v}$.
Assuming $L$ to be the number of hidden layers and $\mathbf{x}^{(0)}=\mathbf{x}\in\mathbb{R}^{n_0}$ the input vector, the NN for $\mathbf{f}(\mathbf{x})$ can be formally expressed as:
\begin{equation}\label{eq:nnout}
\hat{\mathbf{f}}(\mathbf{x}) = \boldsymbol{\Sigma}^{(L+1)} \circ \boldsymbol{\Sigma}^{(L)} \circ \cdots \circ \boldsymbol{\Sigma}^{(1)}(\mathbf{x}).
\end{equation}

The quality of the NN depends on the choice of the weights and biases, which are tuned by minimizing an appropriate loss function, typically defined as the mean squared error of $\hat{\mathbf{f}}$ over the training data points in regression tasks. For classification purposes, the main loss function is the cross-entropy loss. 

The minimization is usually performed by a Stochastic Gradient Descent (SGD) method which iteratively computes the local gradient of the loss function and moves in its descending direction looking for the loss minimum. 
At each SGD iteration (epoch), the method splits the training dataset into small shuffled subsets (mini batches), computes the gradient for each batch, and consequently changes the NN parameters (weights and biases) to move close to the global minimum.  

In our application, the NN vectorial output $\hat{\mathbf{f}}$ is $\mathbf{y}(t)=\{A_a,\delta_S\}$. Note that, the surrogate model was designed to predict $A_a$ and $\delta_S$ rather than $M_0$ directly~\eqref{eq:M0}, as this choice ensures the injectivity of the function being approximated. The NN input is the vector of the parameters $\mathbf{p}=\{\tau_0,\phi,M_2\}$ and the time instant $t$, i.e., $\mathbf{x}^{(0)}=\{\mathbf{p},t\}\in\mathbb{R}^4$. 
By treating time as an explicit input, our model can independently predict outputs for each time step within a fixed time interval. It is beyond the scope of this work to build a time-series forecast of the surrogate model, for which iterative auto-regression models may be employed. An example of such models is provided by long short-term memory (LSTM) networks \cite{HocSch97}, which could offer advantages in scenarios where long-term predictions or auto-regressive dependencies are required.

\begin{figure*}
\centerline{\includegraphics[width=.5\textwidth]{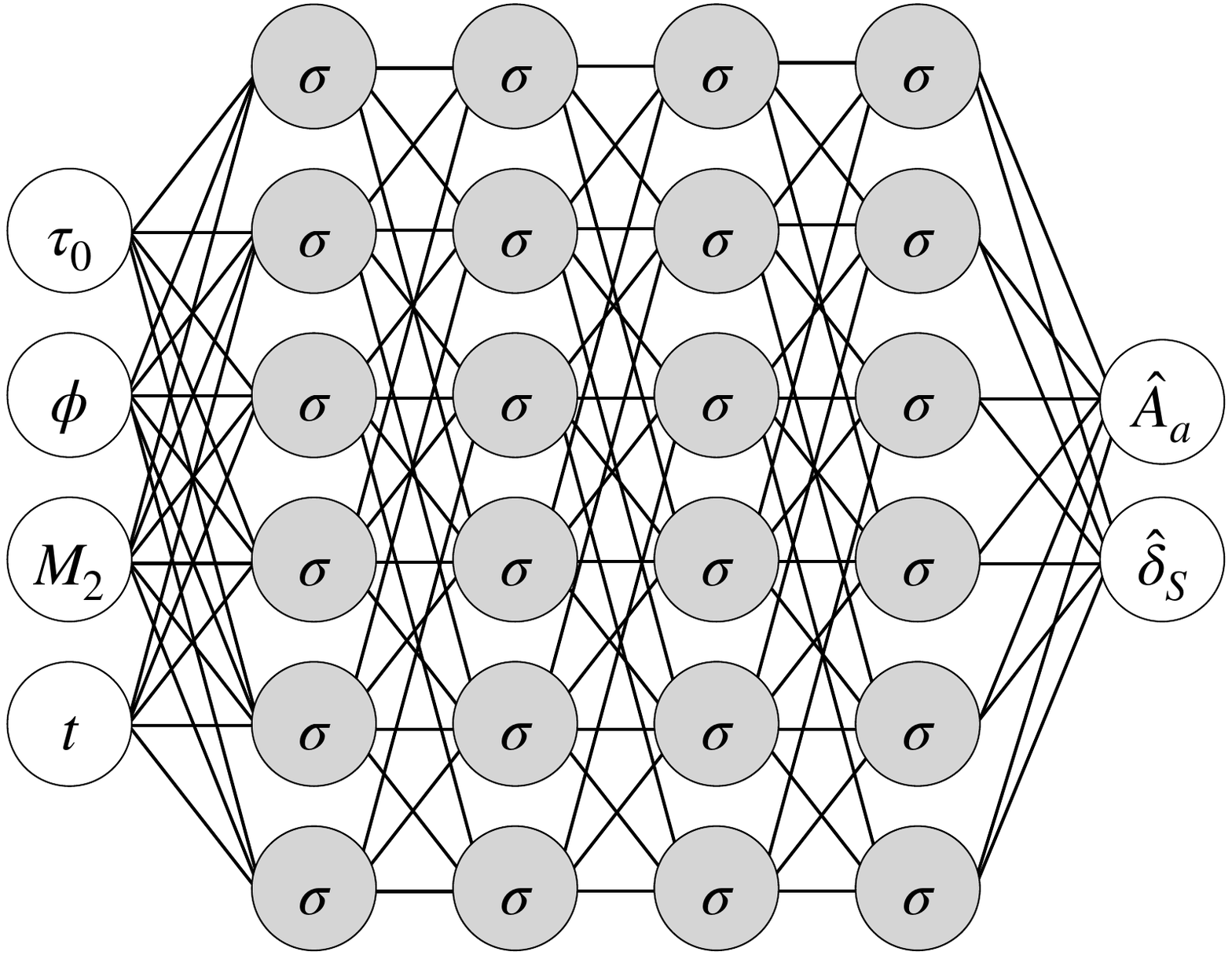}}
\caption{Sketch of a SurMoDeL NN with $L=4$ hidden layers, $n_0=4$, $n_5=2$, and $n_l=6$ for $l=1,\dots,4$.\label{fig:NN}}
\end{figure*}

The set of applied loads $\mathbf{F}$ could be also considered within the input entries, but for the problem at hand we will consider building a surrogate model for some fixed geometry and forcing conditions.
Therefore, the SurMoDeL NN must have a four-dimensional input and a two-dimensional output, so $n_0=4$ and $n_{L+1}=2$. Figure~\ref{fig:NN} shows a sketch of a NN satisfying these requirements.

The SurMoDeL NN is trained by minimizing with a SGD-based method the loss function:
\begin{equation}
    \mathcal{L}(\hat{\mathbf{f}}) = \sum_{i=1}^{N_t} \sum_{j=1}^{N_p} \left\| \hat{\mathbf{f}}(\mathbf{p}_j,t_i) - \mathbf{y}_j(t_i) \right\|_2^2,
    \label{eq:loss}
\end{equation}
for $N_p$ realizations of the uncertain parameter vector $\mathbf{p}$ and $N_t$ time instants. The total amount of data used for the training is therefore $N_d=N_p\times N_t$, which implies running $N_p$ simulations of the full forward model for $N_t$ time-steps each, spanning the time domain $[0,t_{\max}]$.

The architecture of the SurMoDeL NN depends on a number of hyperparameter values, among which the most influential are usually the number $L$ of hidden layers, the number $n_l$ of neurons per layer, and the type of activation function $\sigma^{(l)}$. 
In order to determine the most effective architecture, a sensitivity analysis can be performed to find the best hyperparameter set. 
Under the hypothesis to have the same activation function and number of neurons in each hidden layers, the hyperparameter space is defined as $\mathcal{H}=\mathcal{H}_{L}\times\mathcal{H}_{n_l}\times\mathcal{H}_{\sigma^{(l)}}\times\mathcal{H}_{\sigma^{(L+1)}}$, where $\mathcal{H}_L$, $\mathcal{H}_{n_l}$, $\mathcal{H}_{\sigma^{(l)}}$, and $\mathcal{H}_{\sigma^{(L+1)}}$ are the search spaces
for the number of hidden layers $L$, the number of neurons per layer $n_l$, the type of activation of the hidden layers $\sigma^{(l)}$, and the type of activation function of the output layer $\sigma^{(L+1)}$, respectively. 

\section{Full synthetic model set-up}
\label{sec:setup}

\begin{figure*}
\centerline{\includegraphics[width=\textwidth]{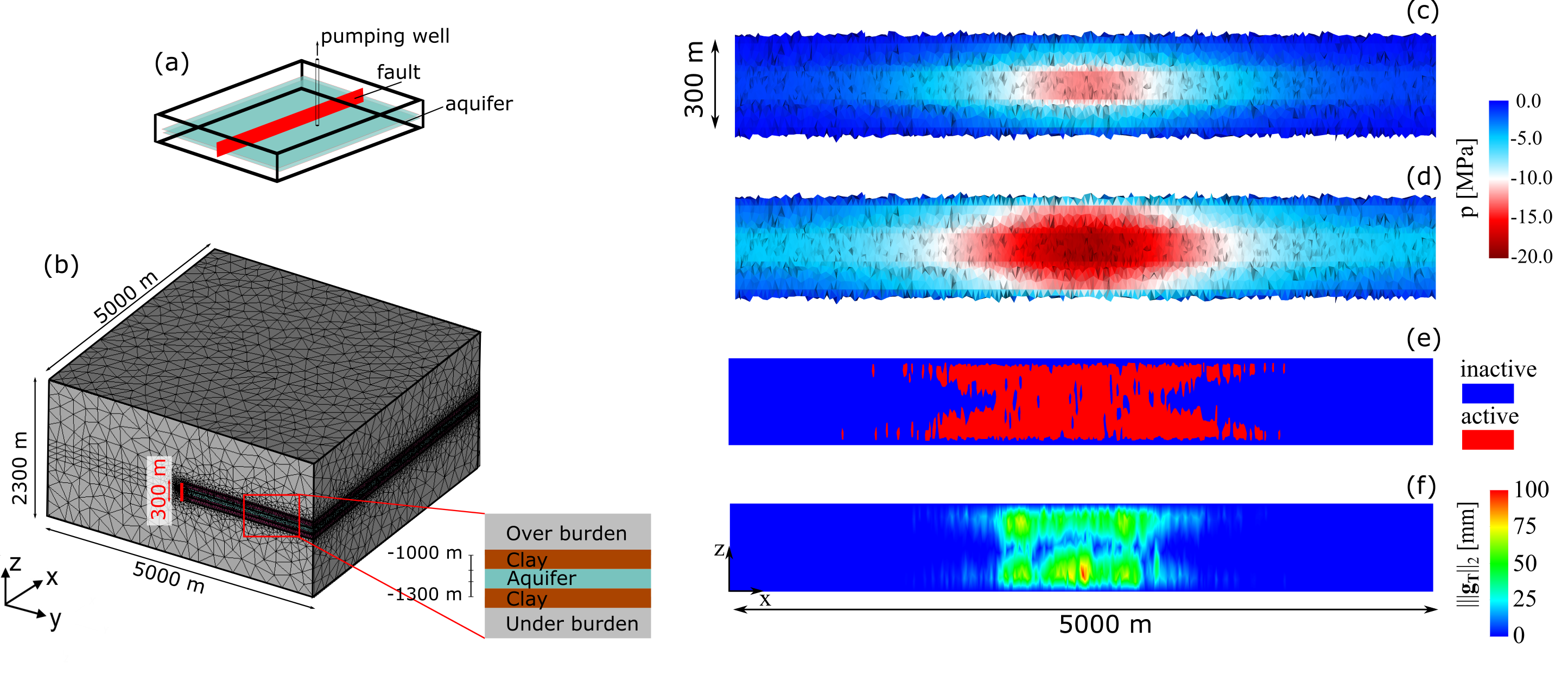}}
\caption{(a-b) Model domain and computational grid used in the full forward simulation. The pumping well produces water from a confined aquifer between -1.100m and -1.200 m. (c-d) Pore pressure distribution within the vertical fault plane at $t_5$ and $t_{10}$. (e-f) Distribution
of sliding (active) and non-sliding (inactive) triangular elements generated by the triangulation over $\Gamma_f$ with the associated sliding values, $\left\| \mathbf{g}_T \right\|_2$, within the vertical fault at $t_{10}$. These outcomes are obtained by running a full forward simulation with the parameter set $\mathbf{p}=\{0,20,0.4286\}$~\cite{Zoc_etal19}. 
\label{fig:modelsetup}}
\end{figure*}

The 3D synthetic case shown in Figure~\ref{fig:modelsetup}a-b and taken from \cite{Zoc_etal19} is used to test and validate the proposed approach, i.e.,
to build the surrogate solution and invert the parameter set by seismic data assimilation. 
It represents an aquifer system cut by a single fault subjected to groundwater abstraction.
As a first step in the workflow, a fluid-dynamical model solving numerically equation \eqref{eq:flow} for a single-phase system is run within the 3D faulted domain to obtain the pore-pressure distribution (Figure~\ref{fig:modelsetup}c-d). The pore pressure outcome has been used as the external source of strength in the full forward geomechanical model~\eqref{eq:forces}. 
The domain extends for 5 km along the $x-$ and $y-$ directions, down to a total depth of 2.300 m.
A discharge of approximately 864 m$^3$/day is constantly pumped from a producing well located 300-m far from the fault in a symmetric position relative to the
$x-$axis over the entire simulation interval of 10 years. Zero-flux boundary conditions are imposed at the domain boundaries. The hydraulic conductivity is equal to 10$^{-7}$ m/s in the aquifer and 10$^{-10}$ m/s in the clay layer within the underburden, siderburden, and overburden. Poisson ratio and Young's modulus are uniform and constant, equal to 0.30 and 1.0 GPa, respectively. 

The forward model \eqref{eq:KKT} is solved by using a tetrahedral discretization of the domain $\Omega$, with the traction over the fault surface $\Gamma_f$ defined by a piecewise constant interpolation carried out on the dual grid generated by the triangulation over $\Gamma_f$.
The overall grid used in the full forward simulation consists of 125,411 nodes and 763,269 elements, with
3,786 triangles discretizing the fault surface. The mesh is particularly refined in the surroundings of the fault and the reservoir (Figure~\ref{fig:modelsetup}b). Boundary conditions are prescribed such that
no displacements are allowed on the bottom boundary
and horizontal displacements are prevented on the lateral
boundaries. The top of the domain is modeled as a traction-free boundary representing the ground surface. 
The simulation spans a temporal interval of ten-time units, hence $N_t=10$ and $t_i=i$, $i=1,\dots,10$. 

A full model run with a deterministic set of parameters $\mathbf{p}=\{0,20,0.4286\}$ is presented in Figure~\ref{fig:modelsetup}e-f, depicting active/inactive triangular elements at $t_{10}$ and the corresponding values $\left\| \mathbf{g}_T \right\|_2$. The majority of the fault sliding occurs in the central portion of the fault, where the pore pressure reaches the maximum values. In particular, Figure~\ref{fig:modelsetup}f shows the distribution of $\left\| \mathbf{g}_T \right\|_2$, highlighting their primary orientation along the $z-$axis with the highest values located at the top and bottom of the fault. These regions correspond to the portions of the fault experiencing the maximum vertical displacements due to aquifer compaction. 

\section{SurMoDeL training and validation}
\label{sec:training}

The goal of the SurMoDeL training is to enable it to learn the complex relationships between the selected inputs and outputs within the geomechanical system. 
The training dataset for our DL-based surrogate model is constructed from points obtained by spanning the parameter space defined in Section~\ref{sec:paramspace} and running the corresponding full forward geomechanical model (Section~\ref{sec:numsim}).

The selected parameter space is the cube $\Psi\equiv\mathcal{D}_{\tau_0}\times\mathcal{D}_{\phi}\times\mathcal{D}_{M_2}$, where: 
\begin{equation}\label{eq:paramdom}
    \mathcal{D}_{\tau_0}=[0,0.2]\text{MPa}, \qquad \mathcal{D}_{\phi}=[20,40]^{\circ}, \qquad
    \mathcal{D}_{M_2}=[0.4286,1.0].
\end{equation}
The cube $\Psi$ is spanned by selecting 5 points per direction, corresponding to the projection over each interval of the Gauss quadrature points in $[-1,1]$. For each one of the $N_p=5^3=125$ combinations, the full forward model is run, getting $A_a$ and $\delta_S$ at each time instant $t_i$. The overall size of the training dataset is therefore $N_d=N_p\times N_t=1250$.  

The $20\%$ of $N_d$ is used as test set, while the remaining part is split into training set and validation set in a ratio of $9$ to $1$. To train the SurMoDeL the maximum number of epochs is set to $10^4$ with an early-stopping condition that ends the training when the loss value over the validation data points (validation loss) has ceased improving for $200$ epochs. Once trained, the model has been evaluated on the whole data set. 
The selected SurMoDeL architecture is defined by running
a random search algorithm over the hyperparameter space $\mathcal{H}$, with $\mathcal{H}_L=\{4,8,12,\dots,40\}$, $\mathcal{H}_{n_l}=\{4,12,20,\dots,100\}$, $\mathcal{H}_{\sigma^{(l)}}=\{\text{ReLU},\tanh,\text{softplus}\}$, $\mathcal{H}_{\sigma^{(L+1)}}=\{\text{ReLU},\text{softplus}\}$. Note that $\mathcal{H}_{\sigma^{(L+1)}}$ does not include the $\tanh$ activation, since the quantities of interest $A_a$ and $\delta_S$ assume only positive values.
The best model, i.e., the one that 
generalizes better, has been identified to be a NN with $L=8$, $n_l=76$, $\sigma^{(l)}=\text{ReLU}$, and $\sigma^{(L+1)}=\text{softplus}$.

\begin{figure*}
\centerline{\includegraphics[width=.5\textwidth]{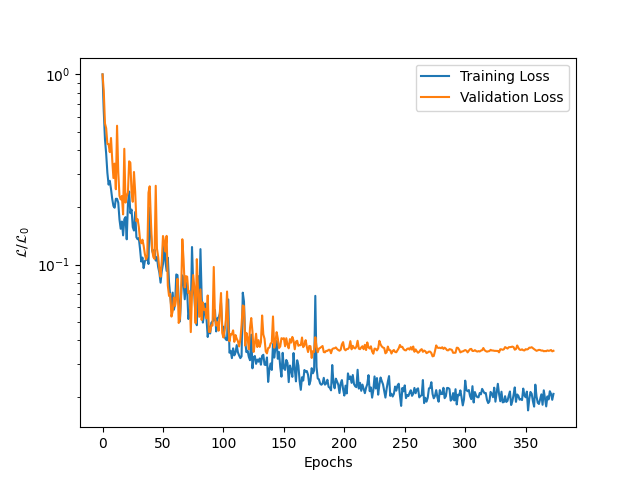}}
\caption{Evolution of the relative loss function during the training of the SurMoDeL.
\label{fig:surmodel_loss}}
\end{figure*}

\begin{figure*}
\centerline{\includegraphics[width=.9\textwidth]{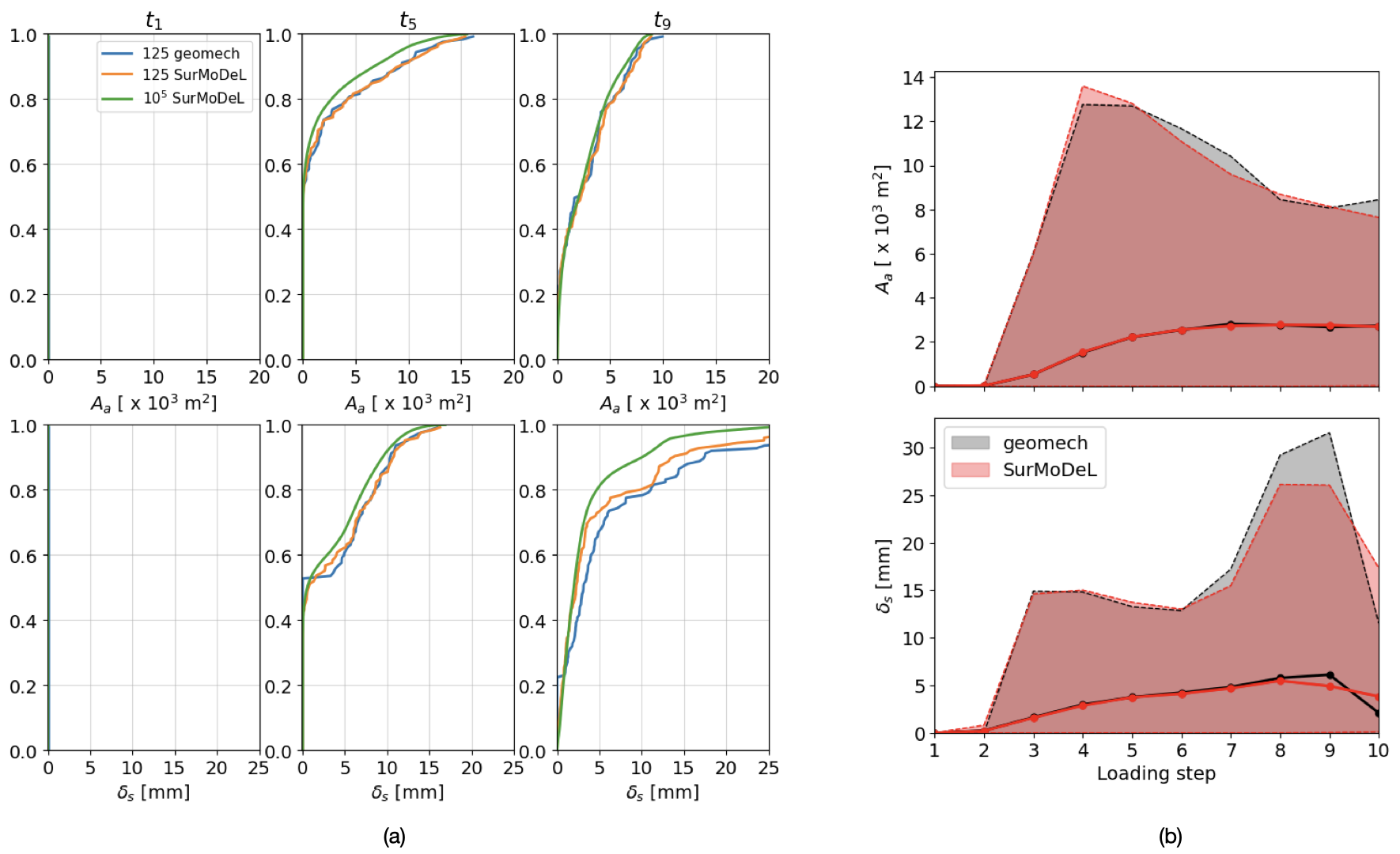}}
\caption{SurMoDeL training results. (a) Cumulative distribution functions of $A_a$ (top row) and $\delta_S$ (bottom row) at different time steps ($t_1$, $t_5$, and $t_9$). The blue lines represent results from the geomechanical model using the $N_p=125$ parameter combinations, while the orange lines depict the outcomes from the SurMoDeL using the same inputs. The green lines show the cumulative distributions from $10^5$ SurMoDeL evaluations on MC realizations. 
(b) Median values (solid lines) and the $2.5\%$ and $97.5\%$ quantiles (dashed lines) for $A_a$ (top) and $\delta_S$ (bottom) obtained using the full forward model (grey) and SurMoDeL (red). \label{fig:res_train}}
\end{figure*}

\begin{center}
\begin{table*}
\caption{Architecture and training conditions of the neural networks. \label{tab:training}}
\begin{tabular*}{\textwidth}{@{\extracolsep{\fill}}cccc@{}}
\toprule
\multicolumn{4}{@{}c}{\textbf{NN Architecture}}
\\
\midrule
$\boldsymbol{\#}$\textbf{ Layers}  & $\boldsymbol{\#}$\textbf{ Neurons per Layer}  & \textbf{Hidden activation} & \textbf{Output Activation} \\
$L=8$  & $n_l=76$  & $\sigma^{(l)}=\text{ReLU}$ & $\sigma^{(L+1)}=\text{softplus}$ \\
\toprule
\multicolumn{4}{@{}c}{\textbf{Dataset}}
\\
\midrule
\textbf{Size}  & \textbf{Training}  & \textbf{Validation} & \textbf{Test} \\
$N_d=1250$  & $72\%$  & $8\%$ & $20\%$ \\
\toprule
\multicolumn{4}{@{}c}{\textbf{Training}}
\\
\midrule
\textbf{Optimizer}  & \textbf{Learning Rate}  & \textbf{Batch Size} & \textbf{Early Stopping} \\
Adam  & $0.001$  & $32$ & Patience of $200$ epochs \\
\bottomrule
\end{tabular*}
\begin{tablenotes}
\item The training is performed in the TensorFlow framework on a machine with two Intel(R) Xeon(R) E5-2680 v2 CPUs @ 2.80GHz and 256GB of RAM.
\end{tablenotes}
\end{table*}
\end{center}
The evolution of the loss~\eqref{eq:loss}, under the training conditions summarized in Table~\ref{tab:training}, is displayed in Figure~\ref{fig:surmodel_loss}, while Figure~\ref{fig:res_train} shows the qualitative results of the training.
Figure~\ref{fig:res_train}a provides the cumulative distribution functions of $A_a$ and $\delta_S$
at different time instants. The outcome obtained by using the full forward geomechanical model on the $N_p=125$ simulations (blue line) is compared to the trained SurMoDeL results for the same set of simulations (orange line), providing a good match. SurMoDeL is also applied on $10^5$ Monte Carlo (MC) samples in the parameter space in~\eqref{eq:paramdom}, getting a uniform distribution (green line) at almost zero cost, since the inference cost of NNs is negligible with respect to one full geomechanical simulation. 
Figure~\ref{fig:res_train}b shows the median of $A_a$ and $\delta_S$ (solid lines) for the ensemble of $N_d=125$ realizations obtained with full forward model and SurMoDeL, along with the $2.5\%$ and $97.5\%$ quantiles (dashed lines) at each time step.
The outcome achieved with the two approaches is very consistent, providing a first validation of the proposed surrogate model. 
From a physical point of view, Figure~\ref{fig:res_train} tells that the fault remains inactive until $t_2$ for any parameter combination. Then, the size of the active area $A_a$ starts increasing in time as the pressure change propagates toward the vertical fault, and the same for the average slip $\delta_S$. After $t_5$, the 97.5\% quantile line for $A_a$ decreases, showing that at this point the fault can also open for some parameter combination. In fact, when the fault opens a portion of $\Gamma_f^{slip}$ becomes $\Gamma_f^{open}$ and this might not be compensated by the portion of $\Gamma_f^{stick}$ that turns into $\Gamma_f^{slip}$.

The SurMoDeL accuracy has been investigated in different training conditions. A subset of the realizations of the uncertain parameter vector of cardinality $N_p=100, 75 \text{ and } 50$ has been randomly selected and then used to build the training data set, thus resulting in a total number of data points $N_d=1000, 750 \text{ or } 500$, respectively.
The SurMoDeL accuracy with these training datasets is compared to the one with the full number of training points in Table~\ref{tab:Ndresults}.
The generalization ability of the proposed SurMoDeL has been analyzed on a new dataset, generated from $N_{MC}=125$ Monte Carlo samples in the parameter domain $\Psi$. 
In particular, Table~\ref{tab:Ndresults} reports: 
\begin{enumerate}
    \item the coefficient of determination: 
    \begin{equation}\label{eq:R2}
        R^2=1-\frac{\sum_{i=1}^{N_t}\sum_{j=1}^{N_{MC}}(y_{\text{ref}}(\mathbf{p}_j,t_i)-\hat{y}(\mathbf{p}_j,t_i))^2}{\sum_{i=1}^{N_t}\sum_{j=1}^{N_{MC}}(y_{\text{ref}}(\mathbf{p}_j,t_i)-\bar{y})^2},
    \end{equation} 
    \item the relative error:
    \begin{equation}\label{eq:E}
        E=\frac{\sum_{i=1}^{N_t}\sum_{j=1}^{N_{MC}}(y_{\text{ref}}(\mathbf{p}_j,t_i)-\hat{y}(\mathbf{p}_j,t_i))^2}{\sum_{i=1}^{N_t}\sum_{j=1}^{N_{MC}}(y_{\text{ref}}(\mathbf{p}_i,t_i))^2},
    \end{equation} 
\end{enumerate}
where $y$ is either the activated area $A_a$ or the average sliding $\delta_S$. Notations $y_{\text{ref}}(\mathbf{p},t)$ and $\hat{y}(\mathbf{p},t)$ refer to the full model and the surrogate model output at the input vector $(\mathbf{p},t)$, respectively, while $\bar{y}=\sum_{i=1}^{N_t}\sum_{j=1}^{N_{MC}}y_{\text{ref}}(\mathbf{p}_j,t_i)/{N_{MC}}$ is the mean of the reference data. 
Table~\ref{tab:Ndresults} also reports the relative error $E$ at each time step for both $A_a$ and $\delta_S$, with the exception of $t_1$ and $t_2$ where the reference solution $y_{\text{ref}}$ is either zero or very close to zero and $E$ is not meaningful. The results in Table~\ref{tab:Ndresults} show that the number of deterministic simulations needed to train the surrogate model can be even reduced to $50$ with a limited loss of accuracy. This is a significant advantage of the proposed approach with respect to other methods used to implement proxy models, such as approximations based on generalized Polynomial Chaos Expansion \cite{Zoc_etal19}, since there is not a minimum number of snapshots needed for the well-posedness of the model assembly. The number of runs with the full forward model used to build the training data set is at the discretion of the modeler, depending on the required proxy model accuracy or the computational cost of the full model. The performance of the surrogate model on the new, previously unseen data set of random simulations demonstrates its generalization ability, despite the small size of the training data set relative to the depth of the selected architecture. Numerical experiments on the unseen data set do not provide evidence of any overfitting issue,
thereby confirming the results of the random search over the NN hyperparameter space that identified the best generalizing architecture.

\begin{center}
\begin{table*}
\caption{Accuracy metrics under different training conditions. \label{tab:Ndresults}}
\begin{tabular*}{\textwidth}{@{\extracolsep\fill}lccccccccc@{}}
\toprule
& & \multicolumn{4}{@{}c}{$A_a$} & \multicolumn{4}{@{}c}{$\delta_S$} \\\cmidrule{3-6}\cmidrule{7-10}
& $t$ & $N_p=125$ & $N_p=100$ & $N_p=75$ & $N_p=50$ & $N_p=125$ & $N_p=100$ & $N_p=75$ & $N_p=50$ \\
\cmidrule{1-2}\cmidrule{3-6}\cmidrule{7-10}
& $t_3$ 
& $5.497\times10^{-3}$  
& $6.120\times10^{-3}$ 
& $3.712\times10^{-3}$  
& $8.768\times10^{-2}$ 
& $4.862\times10^{-2}$  
& $3.315\times10^{-2}$ 
& $4.698\times10^{-2}$  
& $5.337\times10^{-2}$ \\
& $t_4$ 
& $4.577\times10^{-3}$  
& $5.800\times10^{-3}$ 
& $6.703\times10^{-3}$  
& $1.432\times10^{-2}$ 
& $3.873\times10^{-2}$  
& $4.766\times10^{-2}$ 
& $4.100\times10^{-2}$  
& $6.869\times10^{-2}$ \\
& $t_5$ 
& $7.590\times10^{-3}$  
& $7.630\times10^{-3}$ 
& $8.660\times10^{-3}$  
& $1.010\times10^{-2}$ 
& $3.181\times10^{-2}$  
& $3.398\times10^{-2}$ 
& $3.902\times10^{-2}$  
& $4.142\times10^{-2}$ \\
$E$ & $t_6$ 
& $9.722\times10^{-3}$  
& $1.112\times10^{-2}$ 
& $1.528\times10^{-2}$  
& $1.256\times10^{-2}$ 
& $3.947\times10^{-2}$  
& $4.062\times10^{-2}$ 
& $4.071\times10^{-2}$  
& $5.171\times10^{-2}$ \\
& $t_7$ 
& $1.754\times10^{-2}$  
& $1.892\times10^{-2}$ 
& $2.108\times10^{-2}$  
& $2.375\times10^{-2}$ 
& $3.562\times10^{-2}$  
& $6.484\times10^{-2}$ 
& $7.066\times10^{-2}$  
& $9.106\times10^{-2}$ \\
& $t_8$ 
& $2.129\times10^{-2}$  
& $2.325\times10^{-2}$ 
& $2.353\times10^{-2}$  
& $2.901\times10^{-2}$ 
& $5.107\times10^{-2}$  
& $7.239\times10^{-2}$ 
& $1.020\times10^{-1}$  
& $1.071\times10^{-1}$ \\
& $t_9$ 
& $1.765\times10^{-2}$  
& $1.977\times10^{-2}$ 
& $2.033\times10^{-2}$  
& $2.277\times10^{-2}$ 
& $8.450\times10^{-2}$  
& $8.582\times10^{-2}$ 
& $7.391\times10^{-2}$  
& $1.332\times10^{-1}$ \\
& $t_{10}$ 
& $2.867\times10^{-2}$  
& $2.700\times10^{-2}$ 
& $3.254\times10^{-2}$  
& $3.451\times10^{-2}$ 
& $1.363$  
& $11.72$ 
& $13.45$  
& $17.59$ \\
\cmidrule{1-2}\cmidrule{3-6}\cmidrule{7-10}
$E$ & $t_1$-$t_{10}$ 
& $1.602\times10^{-2}$  
& $1.698\times10^{-2}$ 
& $1.898\times10^{-2}$  
& $2.415\times10^{-2}$ 
& $1.143\times10^{-1}$  
& $1.151\times10^{-1}$ 
& $1.286\times10^{-1}$  
& $1.705\times10^{-1}$ \\
$R^2$ & $t_1$-$t_{10}$ 
& $0.979$ 
& $0.977$ 
& $0.975$ 
& $0.968$ 
& $0.846$ 
& $0.845$ 
& $0.827$ 
& $0.770$ \\
\bottomrule
\end{tabular*}
\begin{tablenotes}
\item The coefficient $R^2$ and the mean relative errors $E$ are computed from three runs and different random seeds.
\end{tablenotes}
\end{table*}
\end{center}

\begin{figure*}
\centerline{\includegraphics[width=.9\textwidth]{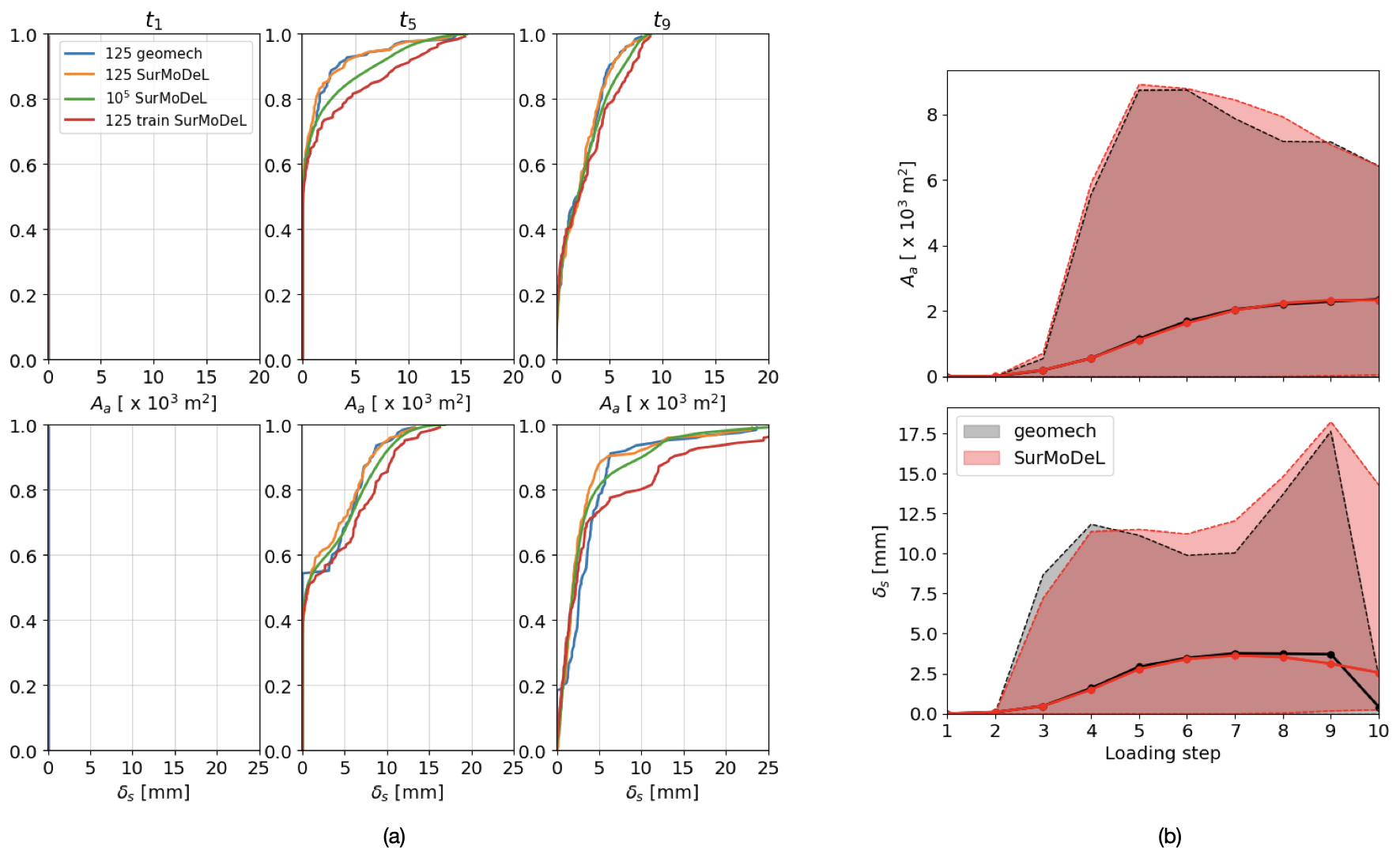}}
\caption{SurMoDeL validation results on 125 random points from the parameter space $\Psi$. (a) Cumulative distribution functions of $A_a$ (top row) and $\delta_S$ (bottom row) at different time steps ($t_1$, $t_5$, and $t_9$). The blue lines represent results from the geomechanical model using the $N_{MC}=125$ MC validation samples, while the orange lines depict the outcomes from the SurMoDeL using the same inputs. The green lines show the cumulative distributions from $10^5$ SurMoDeL evaluations on MC realizations. 
The red line shows for the sake of comparison the SurMoDeL outcome on the training points. (b) Median values (solid lines) and the $2.5\%$ and $97.5\%$ quantiles (dashed lines) for $A_a$ (top) and $\delta_S$ (bottom) obtained using the full forward model (grey) and SurMoDeL (red). \label{fig:res_valid}}
\end{figure*}

The validation of the SurMoDeL training with $N_d=1250$ points is shown in Figure~\ref{fig:res_valid}, which reports the same outcome as Figure~\ref{fig:res_train} computed on a set of 125 random combinations picked from the parameter domain $\Psi$ different from those used for the training. SurMoDeL is able to reproduce almost perfectly the expected behavior of the activated area $A_a$. The approximations in the initial steps for the fault slippage $\delta_S$ are quite accurate as well,
while some challenges appear to arise toward the end of the simulation. In particular, the proposed SurMoDeL is not able to capture satisfactorily the expected decaying trend of $\delta_S$ at $t=10$, as shown also by the large relative errors for $\delta_S$ at $t_{10}$ (Table~\ref{tab:Ndresults}).
As already observed previously, a decrease of $\delta_S$ and $A_a$ with time can be obtained when portions of the fault slip region $\Gamma_f^{slip}$ move to the open region $\Gamma_f^{open}$.  
The proposed surrogate model appears to lack the capability to capture this change in physical behavior. To address this point, it is crucial to develop a proxy model that is able to be aware and classify the different activation modes (slip and open) that can occur. 

\section{Fault activation classification}
\label{sec:faultopen}
\begin{figure*}
\centerline{\includegraphics[width=.9\textwidth]{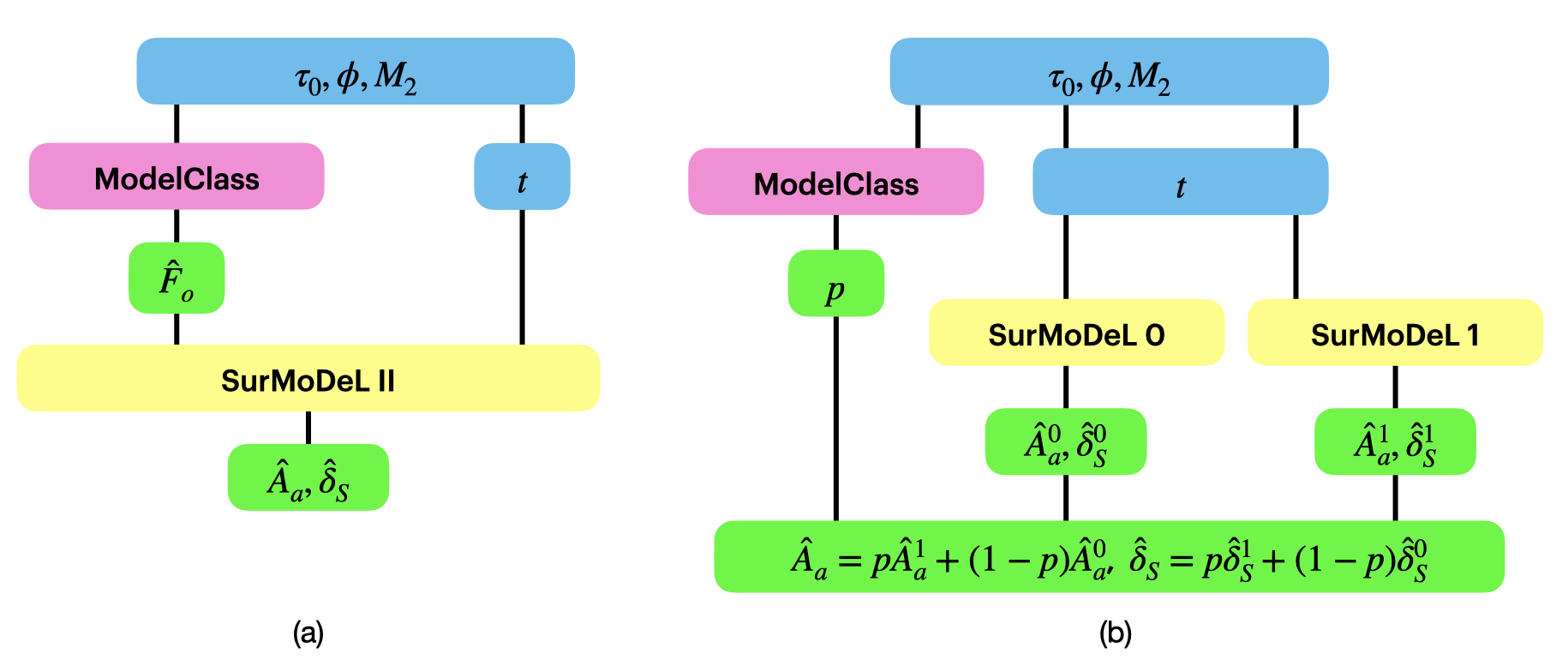}}
\caption{Schematics of the proposed approaches to add physical awareness to the surrogate model:  (a) the logical outcome of the fault activation classification $\hat{F}_o$ is added to the input vector of the surrogate model (SurMoDeL II), and (b) the output of two distinct surrogate models, one trained with no opening occurrence (SurMoDeL 0) and one with opening occurrence (SurMoDeL 1), is combined through the classification prediction $p$. Blue and green colors denote the input and output parameters, respectively, while pink and yellow are the trained NN models. \label{fig:foscheme}}
\end{figure*}

In this section, we discuss algorithmic approaches to enhance the generalization capabilities of the SurMoDeL. The key concept is to introduce some physical awareness in the DL-based surrogate model with the goal of improving the prediction for the last time steps.
To this aim, 
a NN classification model (ModelClass) is trained to foretell when $\Gamma_f^{open}\neq\emptyset$, i.e., when a fault opening occurs. 
Each training point in the parameter space $\Psi$ has been labeled with $1$ if $\Gamma_f^{open}\ne\emptyset$, $0$ otherwise. Hence, the ModelClass takes as input the parameter vector $\mathbf{p}=\{\tau_0,\phi,M_2\}$ and provides as output the probability $p$ of opening occurrence. We set a threshold $\bar{p}$ for such a probability so as to obtain a logical outcome $\hat{F}_o$ for the classifier:
\begin{equation}\label{eq:fo}
    \hat{F}_o=
    \begin{cases}
    1 \quad \text{if } p\ge\bar{p}, \\
    0 \quad \text{if } p<\bar{p}.
    \end{cases}
\end{equation}
Then, the ModelClass prediction is processed in two ways:
\begin{enumerate}
    \item SurMoDeL II: the logical outcome of the fault activation classification $\hat{F}_o$ is added to the surrogate model input vector:
    \begin{equation}
        \hat{\mathbf{y}}=\hat{\mathbf{f}}(\mathbf{p},t,\hat{F}_o);
        \label{eq:surmodelII}
    \end{equation}
    \item SurMoDeL 0 \& SurMoDeL 1: two distinct surrogate models are trained, one with no opening occurrences ($\Gamma_f^{open}=\emptyset$, SurMoDeL 0) and one with opening occurrences ($\Gamma_f^{open}\ne\emptyset$, SurMoDeL 1), with the final output obtained from the linear combinations of the respective outputs $\hat{\mathbf{y}}^0$ and $\hat{\mathbf{y}}^1$ with the probability $p\in[0,1]$ obtained from the fault activation classification:
    \begin{equation}
        \hat{\mathbf{y}} = p \hat{\mathbf{y}}^1 + (1-p) \hat{\mathbf{y}}^0.
    \end{equation}
\end{enumerate}
Figure~\ref{fig:foscheme} illustrates the two different architectures of the new surrogate model that includes the action of ModelClass. Panel (a) represents the single-network approach (SurMoDeL II), where the input parameters $\tau_0$, $\phi$, $M_2$, and time 
$t$ are fed directly into a single neural network to predict the outputs $A_a$ and $\delta_S$. This model also incorporates the output of the ModelClass component, $\hat{F}_0$, in order to account for the awareness of the physical process and refine the prediction capability. Panel (b) depicts the two-network approach, where the classification model first determines the probability $p$ of opening occurrence. Based on this classification, the original training data set is split and two separate surrogate models (SurMoDeL 0 and SurMoDeL 1) are trained to capture distinct behaviors associated to the different physical processes. The final output is a $p$-convex combination of the two NN predictions.

\begin{figure*}
\centerline{\includegraphics[width=.8\textwidth]{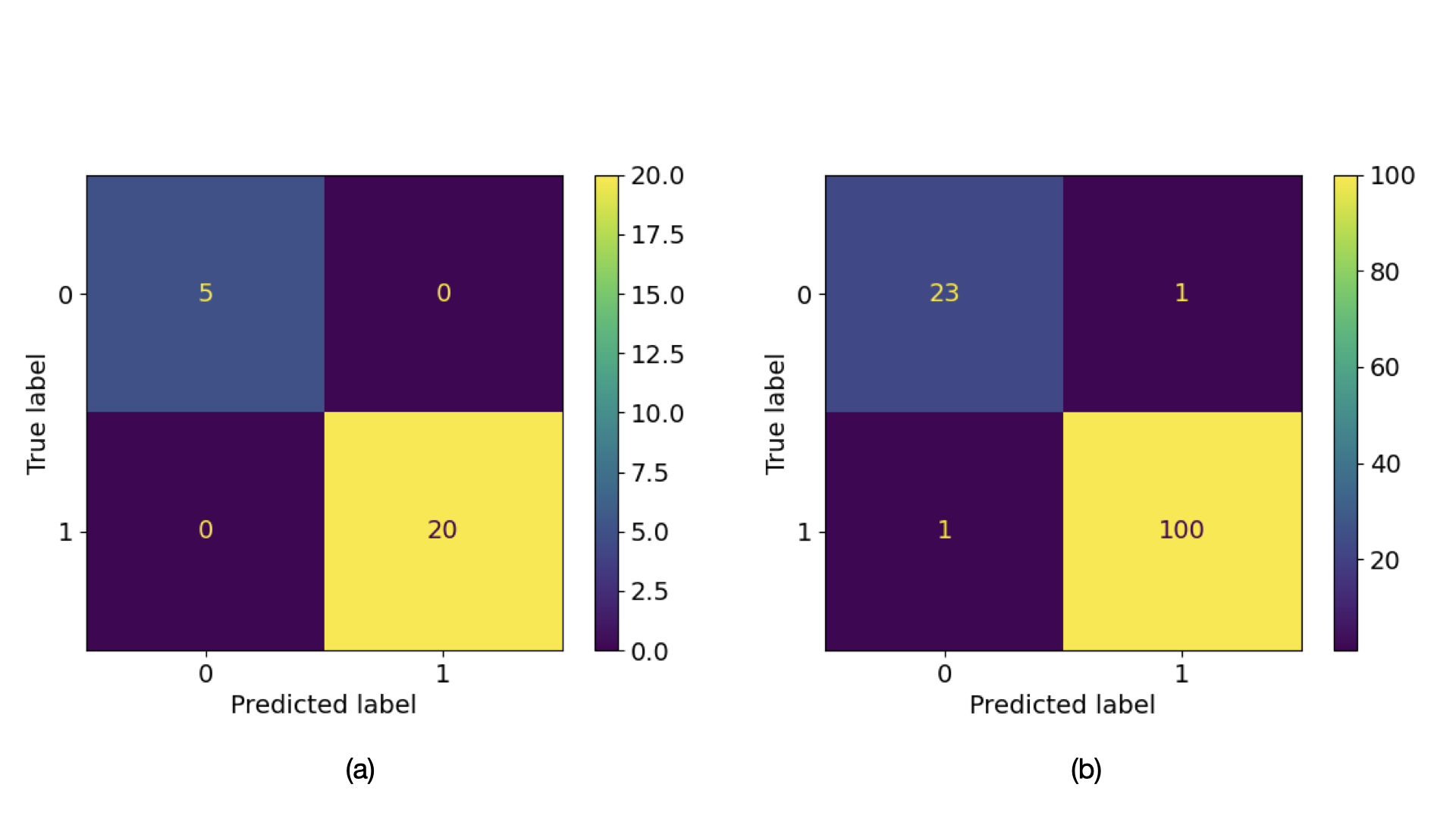}}
\caption{Confusion matrices of ModelClass predictions on (a) the test set, (b) a random set of 125 realizations.\label{fig:cmat}}
\end{figure*}

In our application, we set the threshold $\bar{p}=0.5$.
The NN architecture of ModelClass is derived from a Random Search algorithm over a set of different hyperparameter combinations and finally consists of $16$ layers with $84$ neurons. The hidden activation function is the hyperbolic tangent, while the output activation is the sigmoid function $f(x)=(1+e^{-x})^{-1}$.
ModelClass is trained by using as loss function the binary cross-entropy loss. The accuracy of the classifier is measured by computing how often predictions on the validation dataset match binary labels. We use a batch size equal to $32$ and the training stops before performing $10^4$ epochs if the accuracy metrics does not improve for $200$ epochs.
ModelClass has been finally validated on both the test dataset and the $125$ simulations taken randomly from the parameter space. 
Figure~\ref{fig:cmat} shows the confusion matrices for the predictions on both sets, providing a satisfactory outcome. The diagonal of the matrices represents the number of correctly predicted instances, while the antidiagonal counts the wrong predictions.

\subsection{SurMoDeL II}

\begin{figure*}
\centerline{\includegraphics[width=.9\textwidth]{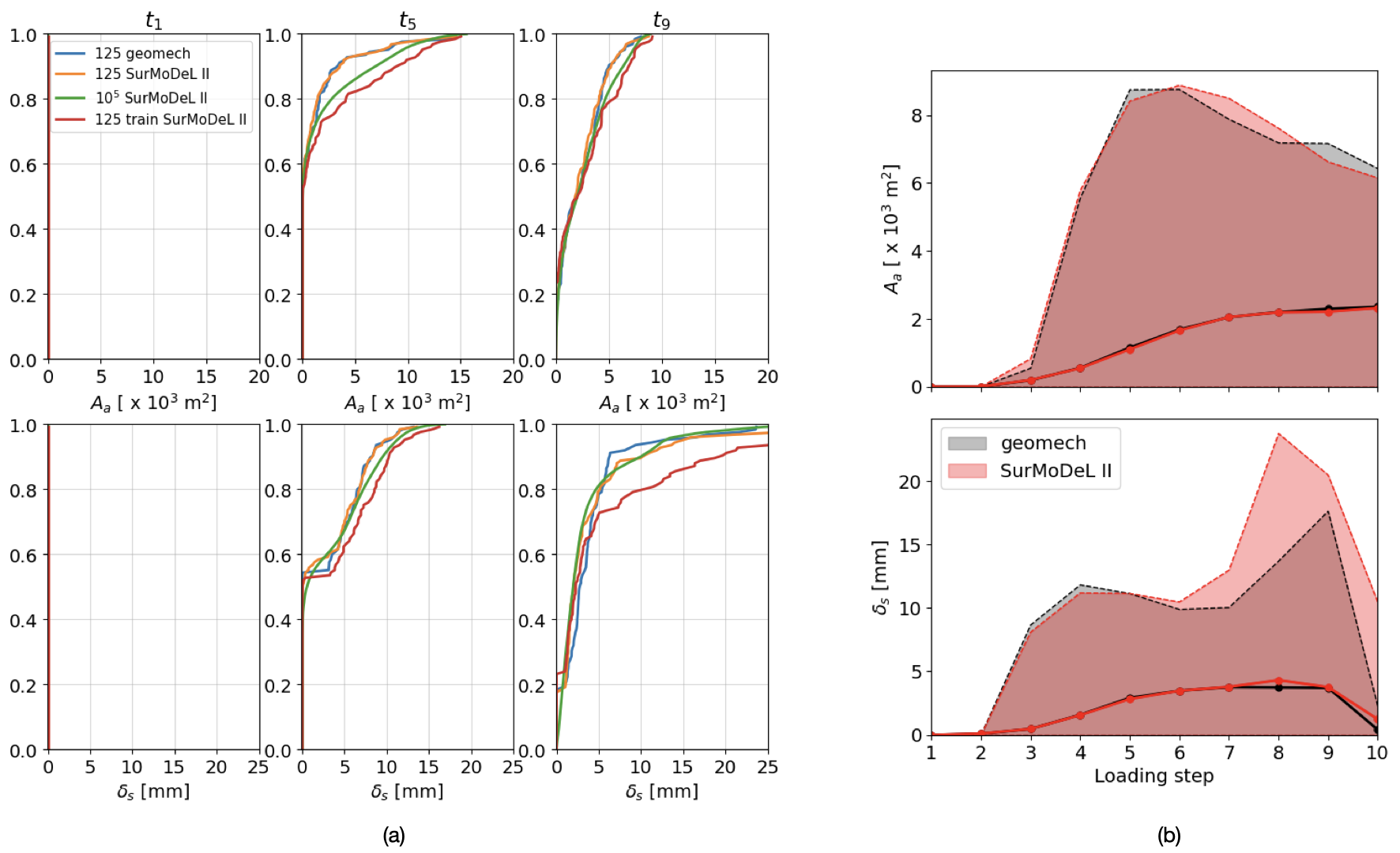}}
\caption{SurMoDeL II validation results on 125 random points from the parameter space $\Psi$. (a) Cumulative distribution functions of $A_a$ (top row) and $\delta_S$ (bottom row) at different time steps ($t_1$, $t_5$, and $t_9$). The blue lines represent results from the geomechanical model using the $N_{MC}=125$ MC validation samples, while the orange lines depict the outcomes from the SurMoDeL II using the same inputs. The green lines show the cumulative distributions from $10^5$ SurMoDeL II evaluations on MC realizations. 
The red line shows for the sake of comparison the SurMoDeL II outcome on the training points. (b) Median values (solid lines) and the $2.5\%$ and $97.5\%$ quantiles (dashed lines) for $A_a$ (top) and $\delta_S$ (bottom) obtained using the full forward model (grey) and SurMoDeL II (red). \label{fig:res_II}}
\end{figure*}

The ModelClass classification is here added to the input of the model. 
The SurMoDeL II design is therefore the same as the surrogate model presented in~\ref{sec:DL}, with the difference in the input vector, which is now of dimension $5$, given that to the parameters $\tau_0$, $\phi$, $M_2$ and the loading time $t$ we add the logical output $\hat{F}_o$.
\begin{figure*}
\centerline{\includegraphics[width=\textwidth]{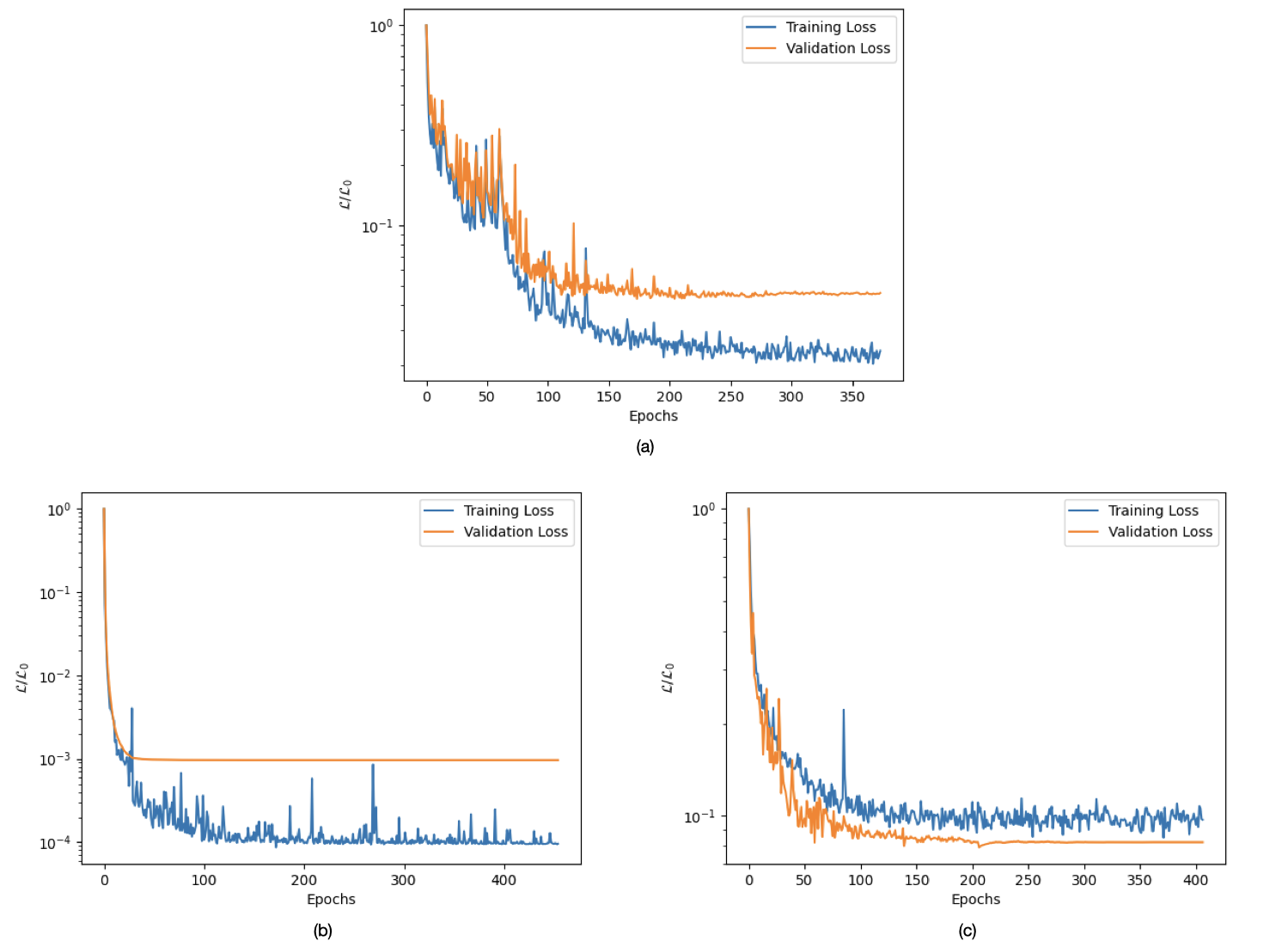}}
\caption{Evolution of the relative loss function during the training of (a) SurMoDeL II, (b) SurMoDeL 0, and (c) SurMoDeL 1.
\label{fig:loss}}
\end{figure*}
The training set has been derived from the one in Section~\ref{sec:training}, by simply evaluating the classifier on each of the training points and adding the corresponding prediction $\hat{F}_o$ to the input vector $\{\mathbf{p},t\}$. The same training conditions as in Section~\ref{sec:training} hold and the relative loss evolution is reported in Figure~\ref{fig:loss}a. The SurMoDeL II validation results on $125$ random points in the parameter space are shown in Figure~\ref{fig:res_II}. In order to evaluate the SurMoDeL II on the validation dataset, we first evaluate ModelClass on $\tau_0,\phi,M_2$, and then use its prediction as extra input of the proxy model. The comparison between the statistics of $\delta_S$ in Figure~\ref{fig:res_valid}b and Figure~\ref{fig:res_II}b shows how the extra input $\hat{F}_o$ impacts on the ability of the proxy model to be aware of the physical fault behavior during the simulation, since now the median and the quantiles start decreasing after timestep $t_8$ as the reference ones. 

\subsection{SurMoDeL 0 \& SurMoDeL 1}

The training data are split on the basis on the fault opening classification, thus generating two datasets of cardinality $N_{d,0}$ and $N_{d,1}$ such that $N_{d,0}+N_{d,1}=N_d=1250$. The first dataset contains all those parameters combinations that do not imply fault opening, the second consists of the remaining triplets $\tau_0,\phi,M_2$ at each timestep $t$. Two distinct surrogate models are created, labeled as SurMoDeL 0 and SurMoDeL 1, which are respectively trained on the first and the second dataset, as described in Section~\ref{sec:training}. Their architecture is defined by the same hyper-parameters of Table~\ref{tab:training}, except for $L=4$, $n_l=38$, and mini-batches of size equal to 16, to account for the lower dimension of the training datasets. The resulting relative losses during the training are provided in Figure~\ref{fig:loss}b and~\ref{fig:loss}c. For any given input vector $\{\mathbf{p},t\}$, SurMoDeL 0 and SurMoDeL 1 predict the related output $\hat{A}_a^0,\hat{\delta}_S^0$ and $\hat{A}_a^1,\hat{\delta}_S^1$, respectively. At the same time, ModelClass provides the probability $p$ associated to the same parameter vector. The outcome of the two surrogate models is finally combined by an affine transformation involving the probability $p$:
\begin{equation}\label{eq:sur01}
    \hat{A}_a=p\hat{A}_a^1+(1-p)\hat{A}_a^0, \quad \hat{\delta}_S=p\hat{\delta}_S^1+(1-p)\hat{\delta}_S^0,
\end{equation}
to obtain the approximations of the activated area $\hat{A}_a$ and the average slippage $\hat{\delta}_S$ of the fault.

\begin{figure*}
\centerline{\includegraphics[width=.9\textwidth]{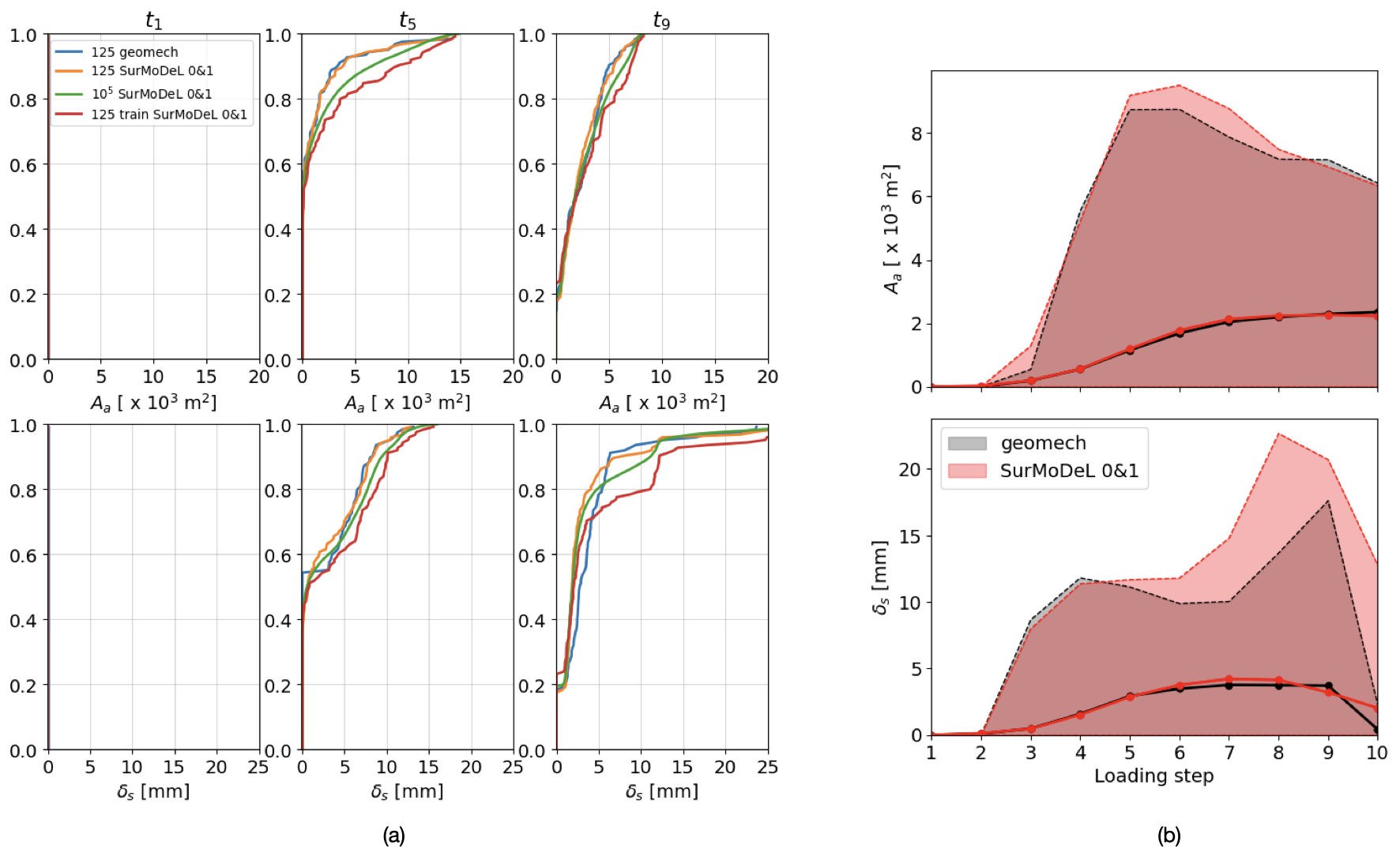}}
\caption{SurMoDeL 0 \& SurMoDeL 1 validation results on 125 random points from the parameter space $\Psi$. (a) Cumulative distribution functions of $A_a$ (top row) and $\delta_S$ (bottom row) at different time steps ($t_1$, $t_5$, and $t_9$). The blue lines represent results from the geomechanical model using the $N_{MC}=125$ MC validation samples, while the orange lines depict the outcomes from the SurMoDeL 0 \& SurMoDeL 1 using the same inputs. The green lines show the cumulative distributions from $10^5$ SurMoDeL 0 \& SurMoDeL 1 evaluations on MC realizations. 
The red line shows for the sake of comparison the SurMoDeL 0 \& SurMoDeL 1 outcome on the training points. (b) Median values (solid lines) and the $2.5\%$ and $97.5\%$ quantiles (dashed lines) for $A_a$ (top) and $\delta_S$ (bottom) obtained using the full forward model (grey) and SurMoDeL 0 \& SurMoDeL 1 (red). \label{fig:res_01}}
\end{figure*}

Figure~\ref{fig:res_01} shows the same results as Figure~\ref{fig:res_II} with SurMoDeL 0 \& SurMoDeL 1. The two diagrams are quite similar, providing evidence that both approaches appear to be effective in adding physical awareness to the proposed surrogate model. A deeper look at Figure~\ref{fig:res_01} shows that SurMoDeL II appears to reproduce the median $\delta_S$ behavior in the final steps slightly better than SurMoDeL 0 \& SurMoDeL 1, but this should not be taken as a general outcome.

\section{Sensitivity Analysis}
\label{sec:SA}

Before moving on to the surrogate model application in seismic data assimilation, we carry out a sensitivity analysis to identify which input parameters are most influential on the uncertainty of the model output. The effects of the input of the model on the output of interest are examined using a variance-based sensitivity analysis~\cite{Sob01}.
The variance of the surrogate model output $\mathbf{y}(t)=\mathbf{f}(\mathbf{p},t)$ can be decomposed as:
\begin{equation}
    \mathrm{Var}(\mathbf{y})=\sum_iV_{\mathbf{p}_i}+\sum_{i<j}V_{\mathbf{p}_i,\mathbf{p}_j}+\cdots+V_{\mathbf{p}_1,\dots,\mathbf{p}_n},
    \label{eq:modelout}
\end{equation}
where
\begin{itemize}
\item $V_{\mathbf{p}_i}=\mathrm{Var}_{\mathbf{p}_i}\left(\mathbb{E}_{{\mathbf{p}}_{\sim i}}(\mathbf{y}\mid \mathbf{p}_{i})\right)$ is the variance with respect to the $i$-th component of the parameter vector $\mathbf{p}_i$ of the expected value of $\mathbf{y}$ taken over all factors but $\mathbf{p}_i$, and represents the contribution to the variance from input $\mathbf{p}_i$;
\item $V_{\mathbf{p}_i,\mathbf{p}_j}=\mathrm{Var}_{\mathbf{p}_i,\mathbf{p}_j}\left(\mathbb{E}_{{\mathbf{p}}_{\sim i},{\mathbf{p}}_{\sim j}}(\mathbf{y}\mid \mathbf{p}_{i},\mathbf{p}_{j})\right)-V_{\mathbf{p}_i}-V_{\mathbf{p}_j}$ is the variance taken over $\mathbf{p}_i$,$\mathbf{p}_{j}$ of the average $\mathbb{E}$ taken over all-but-($\mathbf{p}_i$,$\mathbf{p}_j$) minus their individual input, and captures the contribution from interactions between $\mathbf{p}_i$ and $\mathbf{p}_j$;
\end{itemize}
and so forth. The Sobol indices provide a normalized measure of these contributions~\cite{Sal_etal08}:
\begin{equation}
S_{\mathbf{p}_i} = \frac{V_{\mathbf{p}_i}}{\mathrm{Var}(\mathbf{y})}, \quad S_{\mathbf{p}_i,\mathbf{p}_j} = \frac{V_{\mathbf{p}_i,\mathbf{p}_j}}{\mathrm{Var}(\mathbf{y})}, \quad \dots
\end{equation}

$S_{\mathbf{p}_i}$ are generally called first-order Sobol indices. Note that higher-order Sobol indices do not account for individual contribution of the inputs to the output response. For example, the second-order index $S_{\mathbf{p}_i,\mathbf{p}_j}$ captures the portion of the response of $\mathbf{y}$ to $\mathbf{p}_i$ and $\mathbf{p}_j$ that cannot be expressed as the sum of their separate effects.
In order to total up the contribution of $\mathbf{p}_i$, including all interactions,
the total-effect Sobol index can be computed as: 
\begin{equation}
    S_{T_{\mathbf{p}_i}} = \frac{\mathbb{E}_{{\mathbf{p}}_{\sim i}}\left(\mathrm{Var}_{{\mathbf{p}}_{i}}(\mathbf{y}\mid \mathbf{p}_{i})\right)}{\mathrm{Var}(\mathbf{y})} = 1 - \frac{\mathrm{Var}_{{\mathbf{p}}_{\sim i}}\left(\mathbb{E}_{\mathbf{p}_i}(\mathbf{y} \mid \mathbf{p}_{\sim i})\right)}{\mathrm{Var}(\mathbf{y})},
\end{equation}
where the second equivalence holds since $\sum_iS_{\mathbf{p}_i}+\sum_{i<j}S_{\mathbf{p}_i,\mathbf{p}_j}+\cdots+S_{\mathbf{p}_1,\dots,\mathbf{p}_n}=1$ from~\eqref{eq:modelout}.

Sobol indices are estimated using Monte Carlo-based methods, leveraging sampled input-output data to approximate the necessary conditional variances~\cite{Sob01}. In this study, the SALib module~\cite{HerUsh17,IwaUshHer22} is used to implement the Saltelli extension of the Sobol sequence, a quasi-random low-discrepancy sequence to produce uniform samples of the parameter space~\cite{Sal02}.
The parameter sets are then used to generate the model output using SurMoDeL II, and the sensitivity indices are estimated to allocate output variance to each input.
The Saltelli sampler generates a sample matrix of size $D(2n+2)\times n$, where $D=2^{14}$ and $n=3$ is the number input parameters. The size of $D$ is chosen to ensure the independence of the results from the initialization of the pseudo-random number generator.

Table~\ref{tab:SAresults} presents the Sobol sensitivity indices for the fault activated area ($A_a$) and the average sliding ($\delta_S$) at various loading steps. These indices measure the contributions of the input parameters — $\tau_0$, $\phi$, and $M_2$ — as well as their interactions to the variability of the QoIs. First-order Sobol indices identify the direct influence of the parameter on the output variance, while second order indices reveal the key interactions between two of the parameters.
Regarding the fault activated area, $M_2$ exerts an increasing influence as the loading steps proceed, starting from $S_{M_2}=0.05$ at $t_2$ and growing to $0.70$ by $t_7$ through $t_{10}$. This suggests that $M_2$, hence the stress regime conditions, becomes the primary driver for the fault activated area as the load progresses. The primary contribution of fault cohesion $\tau_0$ increases in the early stages, reaching $0.21$ at $t_{10}$, but remains much lower than $M_2$. The effect of friction angle $\phi$ is consistently small, with $S_{\phi}$ remaining between $0.02$ and $0.04$. The interaction between $M_2$ and $\tau_0$ shows a noticeable contribution early on (with $S_{\tau_0,M_2}=0.19$ at $t_2$), but declines as the loading steps proceed, indicating reduced interaction effects as $M_2$ dominates. The other interaction terms ($S_{\tau_0,\phi}$ and $S_{\phi,M_2}$) are negligible.
$M_2$ also dominates the average sliding, with its influence increasing sharply from $S_{M_2}=0.05$ at $t_2$ to $S_{M_2}=0.86$ at $t_{10}$. This highlights the critical role of $M_2$ in controlling sliding behavior as loading progresses. 
Although significant in $t_3$ ($S_{\tau_0}=0.15$) and $t_4$ ($S_{\tau_0}=0.17$), the influence of $\tau_0$ progressively decreases, reaching near zero from $t_7$.
The contribution of the friction angle is consistently small throughout the loading steps, with $S_{\phi}= 0.02$ or less for most values. The secondary index $S_{\tau_0,M_2}$ assumes notable values in the early steps ($0.21$ at $t_2$ and $0.37$ in $t_3$) but decreases as $M_2$ prevails. Similarly to $A_a$, other interaction terms remain minor.
Note that the sum of the first-order indices is less than one at each loading step, and hence the model is non-additive. 

\begin{center}
\begin{table*}
\caption{Sobol indices of the QoIs. \label{tab:SAresults}}
\begin{tabular*}{\textwidth}{@{\extracolsep\fill}lccccccccccccccccccc@{}}
\toprule
&\multicolumn{9}{@{}c}{$A_a$} & & \multicolumn{9}{@{}c}{$\delta_S$} \\\cmidrule{2-10}\cmidrule{12-20}
 & $t_2$ & $t_3$ & $t_4$ & $t_5$ & $t_6$ & $t_7$ & $t_8$ & $t_9$ & $t_{10}$ & & $t_2$ & $t_3$ & $t_4$ & $t_5$ & $t_6$ & $t_7$ & $t_8$ & $t_9$ & $t_{10}$ \\
\cmidrule{2-10}\cmidrule{12-20}
$S_{\tau_0}$ 
& $0.03$ & $0.09$ & $0.16$ & $0.19$ & $0.19$ & $0.17$ & $0.19$ & $0.20$ & $0.21$ & 
& $0.03$ & $0.15$ & $0.17$ & $0.16$ & $0.12$ & $0.03$ & $0.00$ & $0.01$ & $0.04$ \\
$S_{\phi}$ 
& $0.02$ & $0.03$ & $0.04$ & $0.03$ & $0.03$ & $0.03$ & $0.04$ & $0.04$ & $0.04$ & 
& $0.01$ & $0.02$ & $0.02$ & $0.02$ & $0.01$ & $0.01$ & $0.00$ & $0.00$ & $0.00$ \\
$S_{M_2}$ 
& $0.05$ & $0.16$ & $0.34$ & $0.50$ & $0.62$ & $0.70$ & $0.70$ & $0.69$ & $0.70$ & 
& $0.05$ & $0.30$ & $0.53$ & $0.67$ & $0.75$ & $0.81$ & $0.78$ & $0.81$ & $0.86$ \\
$S_{\tau_0,\phi}$ 
& $0.07$ & $0.08$ & $0.04$ & $0.02$ & $0.01$ & $0.01$ & $0.01$ & $0.01$ & $0.01$ & 
& $0.05$ & $0.02$ & $0.01$ & $0.00$ & $0.00$ & $0.00$ & $0.00$ & $0.00$ & $0.01$ \\
$S_{\tau_0,M_2}$ 
& $0.19$ & $0.34$ & $0.33$ & $0.22$ & $0.13$ & $0.07$ & $0.05$ & $0.05$ & $0.04$ &
& $0.21$ & $0.37$ & $0.18$ & $0.06$ & $0.04$ & $0.10$ & $0.18$ & $0.15$ & $0.07$ \\
$S_{\phi,M_2}$ 
& $0.10$ & $0.09$ & $0.05$ & $0.02$ & $0.01$ & $0.01$ & $0.01$ & $0.00$ & $0.00$ &
& $0.09$ & $0.04$ & $0.02$ & $0.03$ & $0.04$ & $0.02$ & $0.02$ & $0.01$ & $0.01$ \\
\bottomrule
\end{tabular*}
\end{table*}
\end{center}

\begin{figure*}
\centerline{\includegraphics[width=.9\textwidth]{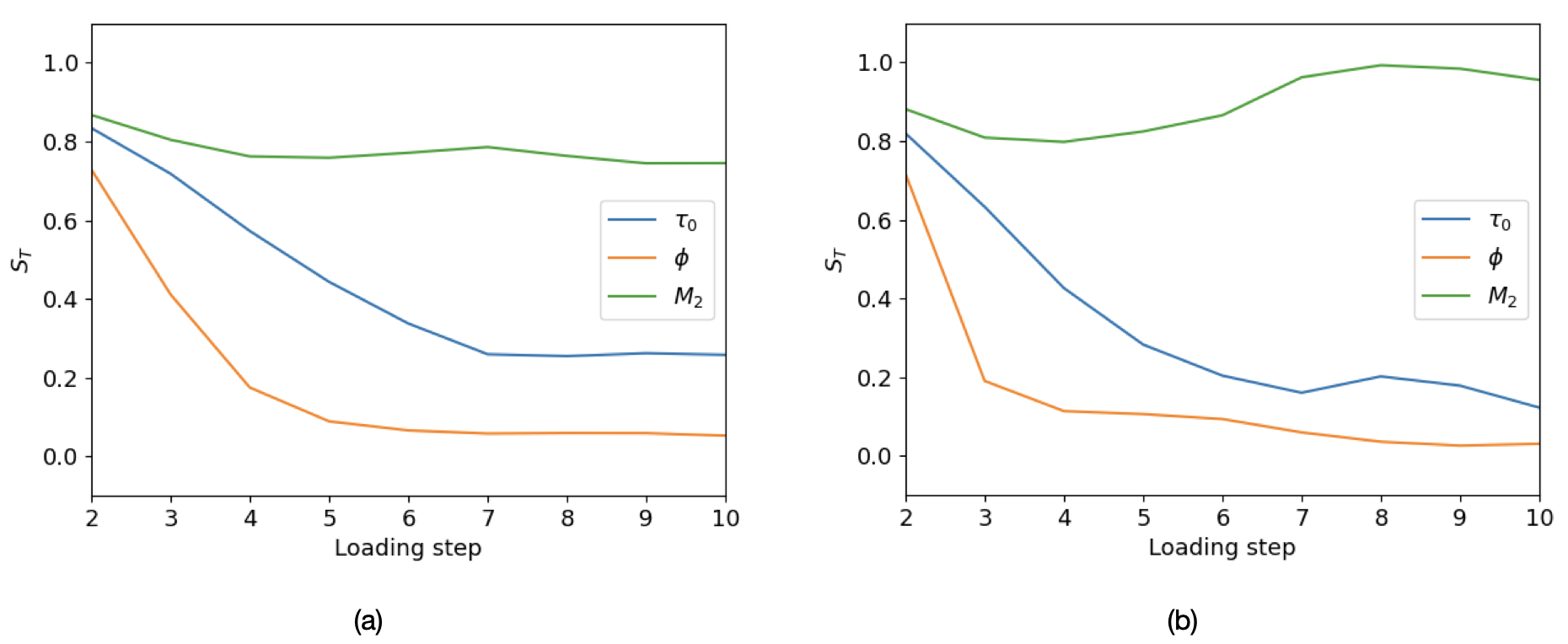}}
\caption{Total effects sensitivity indices ($S_T$) of the input parameters — $\tau_0$, $\phi$, and $M_2$ — on (a) the fault activated area $A_a$ and (b) the average sliding $\delta_S$ for increasing loading steps.\label{fig:SA}}
\end{figure*}

Figure~\ref{fig:SA} shows the total effect $S_T$ of each input parameter in varying loading steps. These indices quantify the extent to which $\tau_0$, $\phi$, or $M_2$, contribute to the variance of the output, accounting for all variance arising from their interactions, of any order, with other input variables. Total effects highlight variables that should not be fixed due to their combined direct and interaction impacts. 
The results show the relative influence of each parameter throughout the progression of the loading, with a consistent decrease in sensitivity for $\tau_0$ and $\phi$, while $M_2$ maintains its prevalence throughout the loading steps. The diminishing sensitivity to $\tau_0$ and $\phi$ suggests that these parameters are less relevant in later stages of the fault response, emphasizing the predominance of $M_2$ in the behavior of the system. 
The fault cohesion $\tau_0$ shows a significant drop in influence, decreasing from $0.7$ in $t_2$ to nearly $0.2$ at $t_{10}$.
Friction angle $\phi$ quickly loses relevance after the initial steps, converging to near-zero total effect by $t_5$. Both QoIs become almost exclusively driven by $M_2$ as loading increases, with its sensitivity index stabilizing at values close to $0.8$ for the fault activated area (Figure~\ref{fig:SA}a) and to $1.0$ for average sliding (Figure~\ref{fig:SA}b). This trend aligns with the results of Table~\ref{tab:SAresults}, confirming the primary role of $M_2$ in determining fault activation and sliding behavior under progressive loading.

\section{Seismic data assimilation}
\label{sec:DA}

The surrogate model is finally used with the aim at solving the inverse problem of estimating $\tau_0$, $\phi$, and $M_2$ from the observation of the seismic moment $M_0$ by a data assimilation approach. 
Data assimilation involves integrating observational data into models to improve their accuracy and reliability. 
Effective data assimilation can bridge the gap between theoretical models and real-world observations, enhancing the model ability to forecast fault activation and the associated seismic risks.

\subsection{Bayesian update}

Bayesian update is a statistical method used to update the probability distribution function (pdf) of model parameters based on new evidence or data. This method is grounded in Bayes' theorem, which describes how to update the probabilities of hypotheses when given evidence.
As it is well-known, Bayes' theorem can be expressed as:
\begin{equation}\label{eq:bayes}
P(\mathbf{p} |\mathbf{q}) = \frac{P(\mathbf{q} | \mathbf{p}) P(\mathbf{p})}{P(\mathbf{q})} 
\end{equation}
where $P(\mathbf{p}|\mathbf{q})$ is the posterior distribution of the parameters $\mathbf{p}$ given the observations $\mathbf{q}$, $P(\mathbf{q}| \mathbf{p})$ is the likelihood of the data given the parameters, $P(\mathbf{p})$ is the prior distribution of the parameters, and $P(\mathbf{q})$ is the marginal likelihood or evidence.

In the context of a Bayesian approach, data assimilation involves the following steps: (i) start with a prior distribution that represents the initial beliefs about the parameters before considering the new data;
(ii) formulate a likelihood function that describes how likely the observed data are, given different values of the parameters; (3) apply Bayes' theorem~\eqref{eq:bayes} to combine the prior distribution and the likelihood function, resulting in the posterior distribution, which reflects the updated beliefs about the parameters after considering the new data. If new data becomes available over time, the posterior distribution from the previous update can serve as the prior distribution for the next update. This process can be repeated as more data are assimilated in time.

By integrating new data with prior information, Bayesian updating provides a rigorous framework for refining model predictions and reducing uncertainty.
The empirical measurements $\mathbf{q}\in\mathbb{R}^{Q}$ are assumed to be noisy versions of the true observable vector $\mathbf{q}_T\in\mathbb{R}^{Q}$:
\begin{equation}
    \mathbf{q} = \mathbf{q}_T + \boldsymbol{\epsilon}
\end{equation}
where $\boldsymbol{\epsilon}\in\mathbb{R}^{Q}$ is the observational error vector, whose components are assumed to be independent and identically distributed with pdf $\pi$. The true values $\mathbf{q}_T=\mathcal{M}\circ\mathcal{S}(\mathbf{p})$ for a given set of loading functions $\mathbf{F}$ represent the observable quantities with the model output $\mathcal{S}(\mathbf{p})$ for the input parameter vector $\mathbf{p}$. If we assume independence between $\boldsymbol{\epsilon}$ and $\mathbf{p}$, the likelihood function is written as:
\begin{equation}
    P(\mathbf{q}|\mathbf{p}) = \prod_{i=1}^{Q} \pi(q_i - \mathcal{M}_i\circ\mathcal{S}(\mathbf{p}))
\end{equation}
This framework allows to combine prior information and empirical data to update the knowledge about the model parameters systematically.

\begin{figure*}
\centerline{\includegraphics[width=.5\textwidth]{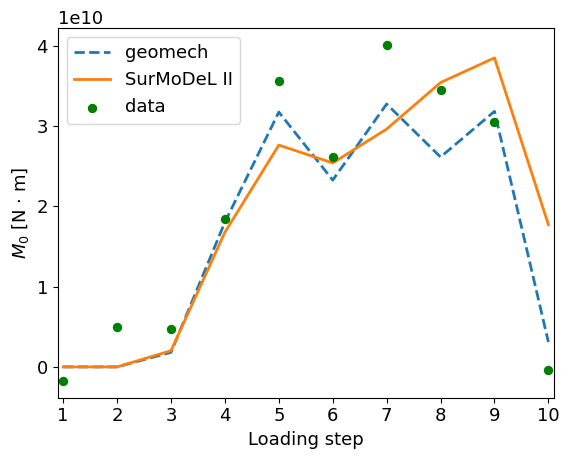}}
\caption{Synthetic test case: seismic moment $M_0$ in time computed using a deterministic run of the full forward model (blue) and the SurMoDeL II approximation (orange). The green dots represent the observation data used in the assimilation process to infer the uncertain parameters $\tau_0$, $\phi$, and $M_2$.\label{fig:M0data}}
\end{figure*}

In this context, we apply the Bayesian inference using the measurement in time of the seismic moment $M_0$ to constrain the model parameters $\mathbf{p}=\{\tau_0,\phi,M_2\}$. 
The observation data $q_i$ used for the Bayesian update are modeled as:
\begin{equation}\label{eq:M0data}
    q_i = M_0(t_i) = M_{0,T}(t_i) + \epsilon_{M_0}(t_i), \quad i=1,\dots,N_t,
\end{equation}
where $M_{0,T}(t)$ is related to the model output through equation~\eqref{eq:M0}, and $\epsilon_{M_0}(t)$ is a Gaussian random noise with standard deviation $\sigma_{\epsilon}=5\cdot10^9$ Nm. 
Since we are testing the proposed approach in the synthetic setting presented in Section~\ref{sec:geom}, we generate the set of "true" observational data $M_{0,T}$ by running the full forward model with the selected "true" parameter set $\tau_0=0.092$ MPa, $\phi=27.1^{\circ}$, $M_2=0.45$. 
They can be inferred from the outcome of the deterministic full model reported as a blue dashed line in Figure~\ref{fig:M0data}. 
The actual $M_{0}(t_i)$ values used in the assimilation process, given by the true reference values disturbed by the noise $\epsilon_{M_0}$, are the green dots in Figure~\ref{fig:M0data}. The orange line corresponds to the outcome predicted by the SurMoDeL II approximations $\hat{A}_a$ and $\hat{\delta}_S$.

The MCMC approach with the Metropolis-Hastings algorithm~\cite{Met_etal53,Has70} are employed as part of the Bayesian inference to estimate the posterior distributions of model parameters, which incorporate both the prior information and the seismic data. This involves iterating the process to sample the posterior distributions effectively. The process begins with an arbitrary initial vector of parameters $\mathbf{p}^{(0)}=\{\tau_0^{(0)},\phi^{(0)},M_2^{(0)}\}$ and a normal pdf is used as the transition kernel in the Metropolis-Hastings algorithm. The prior pdfs for $\tau_0$, $\phi$, and $M_2$ follow uniform distributions: 
\begin{equation}\label{eq:prior}
    \tau_0 \sim U(\mathcal{D}_{\tau_0}), \qquad \phi \sim U(\mathcal{D}_{\phi}), \qquad M_2 \sim U(\mathcal{D}_{M_2}).
\end{equation}
The number of Monte Carlo realizations has been set to $5000$ and SurMoDeL II is employed as a surrogate for evaluating the model outcome, enabling efficient sampling of the posterior distributions without the need for extensive computational resources.

\begin{figure*}
\centerline{\includegraphics[width=.9\textwidth]{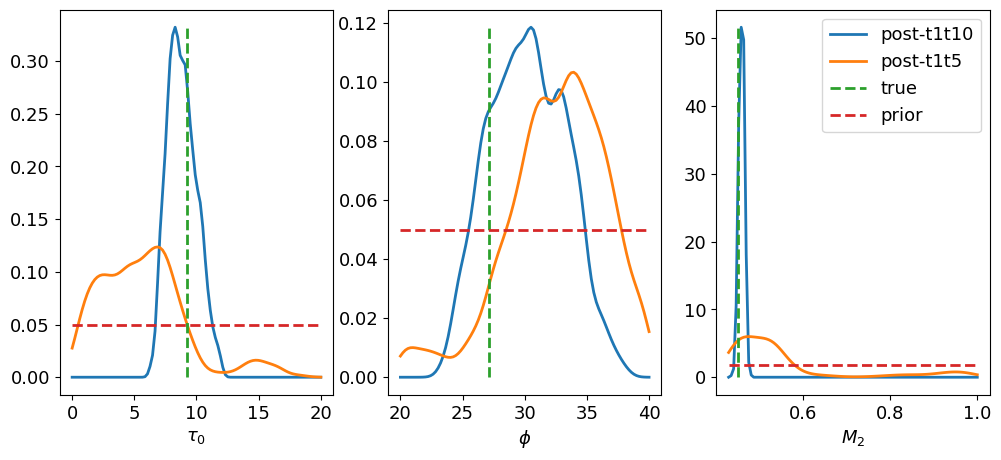}}
\caption{Prior, posterior distributions, and true reference values of the parameters $\tau_0$, $\phi$, and $M_2$ for the two measurement sets considered.\label{fig:post}}
\end{figure*}

Two sets of measurements are used to test the constraint capability:
(i) assimilating $M_0$ data at the first $5$ loading steps ($t_1$-$t_5$),
(ii) assimilating $M_0$ data at all $10$ loading steps ($t_1$-$t_{10}$).
The corresponding posterior distributions, compared to the prior uniform distributions and the true values, are shown in Figure~\ref{fig:post}. 
The initial assimilation in the first half of the time-domain already provides some useful information, even though it appears to be not enough for a satisfactory outcome.
By distinction, the assimilation over the entire dataset is able to provide a very effective constraint on all parameters, and especially so for $\tau_0$ and $M_2$. The obtained results are consistent with the outcome of the sensitivity analysis performed in Section~\ref{sec:SA}, according to which $M_2$ is the most influential parameter, particularly at higher loading steps. 

\section{Conclusions}\label{sec:conclusions}

This work focuses on the development of a theoretical framework to investigate the uncertainties associated with the material parameters in the context of anthropic fault reactivation occurrences. The proposed approach is built on top of the analysis carried out by Zoccarato et al.~\cite{Zoc_etal19}.
A synthetic test case dealing with groundwater extraction from a 1000-m deep confined aquifer bounded by a vertical fault intercepting the sandy-clayey layering system is considered as a reference.

The analysis is aimed at constraining a set of material parameters on the basis of some observable data. In order to solve this inverse problem, a Bayesian inference approach is used.
The main objective of this work is the development of a fast and accurate DL-based surrogate model able to effectively replace the expensive full forward model for the generation of the ensembles of realizations required by the Bayesian assimilation approach.
The following properties of the proposed surrogate model are worth summarizing:
\begin{itemize}
    \item a standard data-driven DL-based approach provides an effective blind alternative to build a surrogate model approximating the functional relationship that connects the output of interest, i.e., the amount of activated area $A_a$ and the average slippage $\delta_S$, to the uncertain input parameters;
    \item the training of the proposed DL-based approach is not very sensitive to the size of the training dataset, providing a similar accuracy also with a relatively small number of points, i.e., a small number of full forward model runs, with respect to other approaches, such as the generalized Polynomial Chaos Expansion~\cite{Zoc_etal19,Zoc_etal20};
    \item introducing some level of awareness of the expected fault physical behavior is very helpful for improving the quality of the DL-based predictions. In the present work, this has been done by connecting the DL-based surrogate model with a prior classifier able to identify the probability of fault opening occurrences as a function of the input parameter set.
\end{itemize}
The parameter space investigated in our analysis is concerned with the estimate of the initial stress regime and the parameters governing the fault failure criterion, while we assume to use as observable data the measurement in time of the seismic moment related to the fault reactivation. A Saltelli global sensitivity analysis underlines the importance of the stress regime in fault activation and sliding behavior, especially under increased loading conditions. The combination of a Bayesian inference carried out by an MCMC approach with the proposed surrogate model turns out to be effective in constraining the models parameters around the "true" values, with a progressive quality increase as the quantity of available data grows. 


\bmsection*{Acknowledgments}

C.M. and M.F. acknowledge the support provided by BIRD2023 and ICEA Department of the University of Padova through the project ``SurMoDeL: Deep Learning Surrogate Models for reservoir characterization'' and their membership of Gruppo Nazionale Calcolo Scientifico - Istituto Nazionale di Alta Matematica (GNCS-INdAM). C.M. is part of the GNCS project CUP\_E53C23001670001 ``Metodi numerici avanzati per la poromeccanica: propriet\`a teoriche ed aspetti computazionali".

\bmsection*{Financial disclosure}

None reported.

\bmsection*{Conflict of interest}

The authors declare no potential conflict of interests.

\bmsection*{Data availability statement}

The data that support the findings of this study are available from the corresponding author upon reasonable request.

\bibliography{biblio}

\begin{thebibliography}{10}
\providecommand \doibase [0]{http://dx.doi.org/}%

\bibitem{Bot_etal19}
Botti M, {Di Pietro} DA, {Le Maître} O, Sochala P. Numerical approximation of
  poroelasticity with random coefficients using Polynomial Chaos and Hybrid
  High-Order methods. {\it Computer Methods in Applied Mechanics and
  Engineering.} 2020\string;361\string:112736.
\newblock \href {\doibase https://doi.org/10.1016/j.cma.2019.112736} {doi:
  https://doi.org/10.1016/j.cma.2019.112736}

\bibitem{BotRos2017}
Bottazzi F, Della~Rossa E. A Functional Data Analysis Approach to Surrogate
  Modeling in Reservoir and Geomechanics Uncertainty Quantification. {\it
  Mathematical Geosciences.} 2017\string;49\string:517--540.
\newblock \href {\doibase 10.1007/s11004-017-9685-y} {doi:
  10.1007/s11004-017-9685-y}

\bibitem{Zoc_etal20}
Zoccarato C, Gazzola L, Ferronato M, Teatini P. Generalized Polynomial Chaos
  Expansion for Fast and Accurate Uncertainty Quantification in Geomechanical
  Modelling. {\it Algorithms.} 2020\string;13(7).
\newblock \href {\doibase 10.3390/a13070156} {doi: 10.3390/a13070156}

\bibitem{Gaz_etal21}
Gazzola L, Ferronato M, Frigo M, et al. A novel methodological approach for
  land subsidence prediction through data assimilation techniques. {\it
  Computational Geosciences.} 2021\string;25(5)\string:1731--1750.
\newblock \href {\doibase 10.1007/s10596-021-10062-1} {doi:
  10.1007/s10596-021-10062-1}

\bibitem{Gaz_etal23}
Gazzola L, Ferronato M, Teatini P, et al. Reducing uncertainty on land
  subsidence modeling prediction by a sequential data-integration approach.
  Application to the Arlua off-shore reservoir in Italy. {\it Geomechanics for
  Energy and the Environment.} 2023\string;33.
\newblock \href {\doibase 10.1016/j.gete.2023.100434} {doi:
  10.1016/j.gete.2023.100434}

\bibitem{Li_etal22}
Li Y, Friedman N, Teatini P, et al. Sensitivity analysis of factors controlling
  earth fissures due to excessive groundwater pumping. {\it Stoch Environ Res
  Risk Assess.} 2022\string;36\string:3911–3928.
\newblock \href {\doibase 10.1007/s00477-022-02237-8} {doi:
  10.1007/s00477-022-02237-8}

\bibitem{CasJhaJua16}
Castiñeira D, Jha B, Juanes R. {Uncertainty Quantification and Inverse
  Modeling of Fault Poromechanics and Induced Seismicity: Application to a
  Synthetic Carbon Capture and Storage (CCS) Problem}.  2016\string;All
  Days\string:ARMA-2016-151.

\bibitem{Zoc_etal19}
Zoccarato C, Ferronato M, Franceschini A, Janna C, Teatini P. Modeling fault
  activation due to fluid production: Bayesian update by seismic data. {\it
  Computational Geosciences.} 2019\string;23(4)\string:705-722.
\newblock \href {\doibase 10.1007/s10596-019-9815-3} {doi:
  10.1007/s10596-019-9815-3}

\bibitem{Verde_2015}
Verde A. Global Sensitivity Analysis of Geomechanical Fractured Reservoir
  Parameters. {\it 49th U.S. Rock Mechanics/Geomechanics Symposium, American
  Rock Mechanics Association.} 2015.

\bibitem{RezNakSid_etal_2020}
Rezaei A, Nakshatrala KB, Siddiqui F, Dindoruk B, Soliman M. A global
  sensitivity analysis and reduced-order models for hydraulically fractured
  horizontal wells. {\it Computational Geosciences.} 2020.
\newblock \href {\doibase 10.1007/s10596-019-09896-7} {doi:
  10.1007/s10596-019-09896-7}

\bibitem{Bla_etal20}
Blaheta R, Béreš M, Domesová S, Horák D. Bayesian inversion for steady flow
  in fractured porous media with contact on fractures and hydro-mechanical
  coupling. {\it Computational Geosciences.} 2020\string;24(5)\string:1911 –
  1932.
\newblock Cited by: 9\href {\doibase 10.1007/s10596-020-09935-8} {doi:
  10.1007/s10596-020-09935-8}

\bibitem{Rut_etal13}
Rutqvist J, Rinaldi AP, Cappa F, Moridis GJ. Modeling of fault reactivation and
  induced seismicity during hydraulic fracturing of shale-gas reservoirs. {\it
  Journal of Petroleum Science and Engineering.} 2013\string;107\string:31-44.
\newblock \href {\doibase https://doi.org/10.1016/j.petrol.2013.04.023} {doi:
  https://doi.org/10.1016/j.petrol.2013.04.023}

\bibitem{Rut_etal16}
Rutqvist J, Rinaldi AP, Cappa F, et al. Fault activation and induced seismicity
  in geological carbon storage – Lessons learned from recent modeling
  studies. {\it Journal of Rock Mechanics and Geotechnical Engineering.}
  2016\string;8(6)\string:789-804.
\newblock \href {\doibase https://doi.org/10.1016/j.jrmge.2016.09.001} {doi:
  https://doi.org/10.1016/j.jrmge.2016.09.001}

\bibitem{Bla_etal22}
Blanco-Martín L, Jahangir E, Rinaldi AP, Rutqvist J. Evaluation of possible
  reactivation of undetected faults during CO2 injection. {\it International
  Journal of Greenhouse Gas Control.} 2022\string;121\string:103794.
\newblock \href {\doibase https://doi.org/10.1016/j.ijggc.2022.103794} {doi:
  https://doi.org/10.1016/j.ijggc.2022.103794}

\bibitem{ChaCheZha10}
Chang H, Chen Y, Zhang D. {Data Assimilation of Coupled Fluid Flow and
  Geomechanics Using the Ensemble Kalman Filter}. {\it SPE Journal.}
  2010\string;15(02)\string:382-394.
\newblock \href {\doibase 10.2118/118963-PA} {doi: 10.2118/118963-PA}

\bibitem{EmeRey12}
Emerick AA, Reynolds AC. {History matching time-lapse seismic data using the
  ensemble Kalman filter with multiple data assimilations}. {\it Computational
  Geosciences.} 2012\string;16(3)\string:639-659.
\newblock \href {\doibase 10.1007/s10596-012-9275-5} {doi:
  10.1007/s10596-012-9275-5}

\bibitem{EmeRey13}
Emerick AA, Reynolds AC. {History-Matching Production and Seismic Data in a
  Real Field Case Using the Ensemble Smoother With Multiple Data Assimilation}.
   2013\string;All Days\string:SPE-163675-MS.
\newblock \href {\doibase 10.2118/163675-MS} {doi: 10.2118/163675-MS}

\bibitem{Luo_etal16}
Luo X, Bhakta T, Jakobsen M, Nædal G. {An Ensemble 4D Seismic History Matching
  Framework with Sparse Representation Based on Wavelet Multiresolution
  Analysis}.  2016\string;Day 1 Wed, April 20, 2016\string:D011S005R002.
\newblock \href {\doibase 10.2118/180025-MS} {doi: 10.2118/180025-MS}

\bibitem{Luo_etal18}
X L, T B, M J, G N. Efficient big data assimilation through sparse
  representation: A 3D benchmark case study in petroleum engineering.. {\it
  PLoS ONE.} 2018\string;13(7).
\newblock \href {\doibase https://doi.org/10.1371/journal.pone.0198586} {doi:
  https://doi.org/10.1371/journal.pone.0198586}

\bibitem{Nej_etal19}
Nejadi S, Kazemi N, Hubbard SM, Gates ID. {A Bayesian Sampling Framework With
  Seismic Priors for Data Assimilation and Uncertainty Quantification}.
  2019\string;Day 1 Wed, April 10, 2019\string:D010S017R009.
\newblock \href {\doibase 10.2118/193844-MS} {doi: 10.2118/193844-MS}

\bibitem{Pan_etal14}
Pan I, Babaei M, Korre A, Durucan S. Artificial Neural Network based surrogate
  modelling for multi- objective optimisation of geological CO2 storage
  operations. {\it Energy Procedia.} 2014\string;63\string:3483-3491.
\newblock 12th International Conference on Greenhouse Gas Control Technologies,
  GHGT-12\href {\doibase https://doi.org/10.1016/j.egypro.2014.11.377} {doi:
  https://doi.org/10.1016/j.egypro.2014.11.377}

\bibitem{She_etal22}
Shen L, Li D, Zha W, Li X, Liu X. Surrogate modeling for porous flow using deep
  neural networks. {\it Journal of Petroleum Science and Engineering.}
  2022\string;213\string:110460.
\newblock \href {\doibase https://doi.org/10.1016/j.petrol.2022.110460} {doi:
  https://doi.org/10.1016/j.petrol.2022.110460}

\bibitem{Qi_etal23}
Qi M, Jang K, Cui C, Moon I. Novel control-aware fault detection approach for
  non-stationary processes via deep learning-based dynamic surrogate modeling.
  {\it Process Safety and Environmental Protection.}
  2023\string;172\string:379-394.
\newblock \href {\doibase https://doi.org/10.1016/j.psep.2023.02.023} {doi:
  https://doi.org/10.1016/j.psep.2023.02.023}

\bibitem{Meg_etal23}
Meguerdijian S, Pawar RJ, Chen B, Gable CW, Miller TA, Jha B. Physics-informed
  machine learning for fault-leakage reduced-order modeling. {\it International
  Journal of Greenhouse Gas Control.} 2023\string;125\string:103873.
\newblock \href {\doibase https://doi.org/10.1016/j.ijggc.2023.103873} {doi:
  https://doi.org/10.1016/j.ijggc.2023.103873}

\bibitem{Lu_etal23}
Lu H, Salgado LS, Marzouk Y, Juanes R. Uncertainty Quantification of CO 2
  Leakage and Risk Analysis of Induced Seismicity for Large-scale Geological CO
  2 Sequestration. {\it AGU23.} 2023.

\bibitem{KikOde88}
Kikuchi N, Oden JT. {\it {Contact Problems in Elasticity: A Study of
  Variational Inequalities and Finite Element Methods}}.
\newblock Philadelphia, PA, USA: SIAM, 1988

\bibitem{Lau03}
Laursen TA. {\it {Computational Contact and Impact Mechanics: Fundamentals of
  Modeling Interfacial Phenomena in Nonlinear Finite Element Analysis}}.
\newblock Springer-Verlag Berlin Heidelberg, 2003

\bibitem{Wri06}
Wriggers P. {\it {Computational Contact Mechanics}}.
\newblock Springer-Verlag Berlin Heidelberg.
\newblock 2nd~ed., 2006

\bibitem{SimHug98}
Simo JC, Hughes TJR. {\it {Computational Inelasticity}}.
\newblock Springer-Verlag New York, 1998

\bibitem{Ber84}
Bertsekas DP. {\it {Constrained Optimization and Lagrange Multiplier Methods}}.
\newblock Academic Press New York, 1984.

\bibitem{HagHueWoh08}
Hager C, H\"{u}eber S, Wohlmuth BI. {A stable energy-conserving approach for
  frictional contact problems based on quadrature formulas}. {\it Int. J.
  Numer. Meth. Eng..} 2008\string;73(2)\string:205--225.
\newblock \href {\doibase 10.1002/nme.2069} {doi: 10.1002/nme.2069}

\bibitem{fraferjantea16}
Franceschini A, Ferronato M, Janna C, Teatini P. {A novel Lagrangian approach
  for the stable numerical simulation of fault and fracture mechanics}. {\it J.
  Comput. Phys..} 2016\string;314\string:503--521.
\newblock \href {\doibase 10.1016/j.jcp.2016.03.032} {doi:
  10.1016/j.jcp.2016.03.032}

\bibitem{fr2020alg}
Franceschini A, Castelletto N, White JA, Tchelepi HA. {Algebraically stabilized
  Lagrange multiplier method for frictional contact mechanics with
  hydraulically active fractures}. {\it Comput. Meth. in Appl. Mech. Eng..}
  2020\string;368\string:113161.
\newblock \href {\doibase 10.1016/j.cma.2020.113161} {doi:
  10.1016/j.cma.2020.113161}

\bibitem{FraFerFriJan22}
Franceschini A, Ferronato M, Frigo M, Janna C. {A reverse augmented constraint
  preconditioner for Lagrange multiplier methods in contact mechanics}. {\it
  Comput. Meth. in Appl. Mech. Eng..} 2022\string;392\string:114632.
\newblock \href {\doibase 10.1016/j.cma.2022.114632} {doi:
  10.1016/j.cma.2022.114632}

\bibitem{Fer_etal13}
Ferronato M, Castelletto N, Gambolati G, Janna C, Teatini P. {II cycle
  compressibility from satellite measurements}. {\it Geotechnique.}
  2013\string;63\string:479-486.
\newblock \href {\doibase 10.1680/geot.11.P.149} {doi: 10.1680/geot.11.P.149}

\bibitem{Zoc_etal16}
Zoccarato C, Ba\`u D, , et al. {Data assimilation of surface displacements to
  improve geomechanical parameters of gas storage reservoirs}. {\it J. Geophys.
  Res. Solid Earth.} 2016\string;121\string:1441-1461.
\newblock \href {\doibase 10.1002/2015JB012090} {doi: 10.1002/2015JB012090}

\bibitem{ZocFerTea18}
Zoccarato C, Ferronato M, Teatini P. {Formation compaction vs land subsidence
  to constrain rock compressibility of hydrocarbon reservoirs}. {\it Geomech.
  Energy Environ..} 2018\string;13\string:14-24.
\newblock \href {\doibase 10.1016/j.gete.2017.12.002} {doi:
  10.1016/j.gete.2017.12.002}

\bibitem{KanAnd75}
Kanamori H, Anderson D. {Theoretical basis of some empirical relations in
  seismology}. {\it Bull. Seismol. Soc. Am..} 1975\string;65\string:1073-1095.

\bibitem{HocSch97}
Hochreiter S, Schmidhuber J. Long Short-Term Memory. {\it Neural Computation.}
  1997\string;9(8)\string:1735-1780.
\newblock \href {\doibase 10.1162/neco.1997.9.8.1735} {doi:
  10.1162/neco.1997.9.8.1735}

\bibitem{Sob01}
Sobol I. Global sensitivity indices for nonlinear mathematical models and their
  Monte Carlo estimates. {\it Mathematics and Computers in Simulation.}
  2001\string;55(1)\string:271-280.
\newblock The Second IMACS Seminar on Monte Carlo Methods\href {\doibase
  https://doi.org/10.1016/S0378-4754(00)00270-6} {doi:
  https://doi.org/10.1016/S0378-4754(00)00270-6}

\bibitem{Sal_etal08}
A S, M R, T A, et al. {\it Global Sensitivity Analysis: The Primer}.
\newblock Chichester (England): Wiley, 2008.

\bibitem{HerUsh17}
Herman J, Usher W. {SALib}: An open-source Python library for Sensitivity
  Analysis. {\it The Journal of Open Source Software.} 2017\string;2(9).
\newblock \href {\doibase 10.21105/joss.00097} {doi: 10.21105/joss.00097}

\bibitem{IwaUshHer22}
Iwanaga T, Usher W, Herman J. Toward SALib 2.0: Advancing the accessibility and
  interpretability of global sensitivity analyses. {\it Socio-Environmental
  Systems Modelling.} 2022\string;4\string:18155.
\newblock \href {\doibase 10.18174/sesmo.18155} {doi: 10.18174/sesmo.18155}

\bibitem{Sal02}
Saltelli A. Making best use of model evaluations to compute sensitivity
  indices. {\it Computer Physics Communications.}
  2002\string;145(2)\string:280-297.
\newblock \href {\doibase https://doi.org/10.1016/S0010-4655(02)00280-1} {doi:
  https://doi.org/10.1016/S0010-4655(02)00280-1}

\bibitem{Met_etal53}
Metropolis N, Rosenbluth AW, Rosenbluth MN, Teller AH, Teller E. {Equation of
  State Calculations by Fast Computing Machines}. {\it The Journal of Chemical
  Physics.} 1953\string;21(6)\string:1087-1092.
\newblock \href {\doibase 10.1063/1.1699114} {doi: 10.1063/1.1699114}

\bibitem{Has70}
Hastings WK. Monte Carlo Sampling Methods Using Markov Chains and Their
  Applications. {\it Biometrika.} 1970\string;57(1)\string:97--109.

\bibitem{Bal_etal24}
Baldan S, Ferronato M, Franceschini A, et al. Unexpected fault activation in
  underground gas storage. Part II: Definition of safe operational bandwidths.
  2024.
\newblock arXiv:2408.01049, math.NA.

\bibitem{Che11}
Chen P. {Full-Wave Seismic Data Assimilation: Theoretical Background and Recent
  Advances}. {\it Pure and Applied Geophysics.}
  2011\string;168(10)\string:1527-1552.
\newblock \href {\doibase 10.1007/s00024-010-0240-8} {doi:
  10.1007/s00024-010-0240-8}

\bibitem{Che_etal20}
Chen M, Abdalla OA, Izady A, {Reza Nikoo} M, Al-Maktoumi A. Development and
  surrogate-based calibration of a CO2 reservoir model. {\it Journal of
  Hydrology.} 2020\string;586\string:124798.
\newblock \href {\doibase https://doi.org/10.1016/j.jhydrol.2020.124798} {doi:
  https://doi.org/10.1016/j.jhydrol.2020.124798}

\bibitem{Fer_etal10}
Ferronato M, Gambolati G, Janna C, Teatini P. Geomechanical issues of
  anthropogenic CO2 sequestration in exploited gas fields. {\it Energy
  Conversion and Management.} 2010\string;51(10)\string:1918-1928.
\newblock \href {\doibase https://doi.org/10.1016/j.enconman.2010.02.024} {doi:
  https://doi.org/10.1016/j.enconman.2010.02.024}

\bibitem{Fra_etal24}
Franceschini A, Zoccarato C, Baldan S, et al. Unexpected fault activation in
  underground gas storage. Part I: Mathematical model and mechanisms. 2023.
\newblock arXiV:2308.02198, math.NA.

\bibitem{Ket_etal16}
Keating EH, Harp DH, Dai Z, Pawar RJ. Reduced order models for assessing CO2
  impacts in shallow unconfined aquifers. {\it International Journal of
  Greenhouse Gas Control.} 2016\string;46\string:187-196.
\newblock \href {\doibase https://doi.org/10.1016/j.ijggc.2016.01.008} {doi:
  https://doi.org/10.1016/j.ijggc.2016.01.008}

\bibitem{OliReyLiu08}
Oliver DS, Reynolds AC, Liu N. {\it Inverse Theory for Petroleum Reservoir
  Characterization and History Matching}.
\newblock Cambridge University Press, 2008.

\end{thebibliography}

\nocite{*}

\end{document}